\documentclass[12pt] {article}
\usepackage{amsmath,amsthm}
\usepackage{amssymb,latexsym}
\usepackage{amscd}
\title{Theory of valuations on manifolds, I. Linear spaces.}
\date{}
\author{ Semyon Alesker \footnote{Partially supported by ISF grant 1369/04.}
\\  { \normalsize Department of Mathematics, Tel Aviv University, Ramat Aviv}
 \\  { \normalsize 69978 Tel Aviv,
Israel }
\\ {\normalsize e-mail: semyon@post.tau.ac.il}}

\def\eps{\varepsilon}

\def\ome{\omega}
\def\Ome{\Omega}
\def\lam{\lambda}

\def\to{\rightarrow}
\def\qed { Q.E.D. }
\def\pt{\partial}

\def\RR{\mathbb{R}}
\def\CC{\mathbb{C}}

\def\NN{\mathbb{N}}
\def\ZZ{\mathbb{Z}}

\def\PP{\mathbb{P}}

\def\ff { \bar {\cal F} }
\swapnumbers
\newtheorem{theorem}{Theorem}[subsection]
\newtheorem{corollary}[theorem]{Corollary}
\newtheorem{lemma}[theorem]{Lemma}
\newtheorem{proposition}[theorem]{Proposition}
\newtheorem{claim}[theorem]{Claim}
\theoremstyle{definition}
\newtheorem{example}[theorem]{Example}
\newtheorem{definition}[theorem]{Definition}
\newtheorem{remark}[theorem]{Remark}
\theoremstyle{proposition-definition}
\newtheorem{proposition-definition}[theorem]{Proposition-Definition}

\def\cf{{\cal F}}

  \def\cc{{\cal C}}
 \def\ce{{\cal E}} \def\cf{{\cal F}}
\def\cg{{\cal G}}  
 \def\ck{{\cal K}} \def\cl{{\cal L}}
  \def\co{{\cal O}}

\def\sv{SV}
\def\svv{SV(V)}

\def\qvv{QV (V)}

\def\pl{\PP_+(V^*)}

\def\mov{|\omega_V|}
\makeatletter

\def\diagram{\m@th\leftwidth=\z@ \rightwidth=\z@ \topheight=\z@
\botheight=\z@ \setbox\@picbox\hbox\bgroup}

\def\enddiagram{\egroup\wd\@picbox\rightwidth\unitlength
\ht\@picbox\topheight\unitlength \dp\@picbox\botheight\unitlength
\hskip\leftwidth\unitlength\box\@picbox}

\def\bfig{\begin{diagram}}
\def\efig{\end{diagram}}
\newcount\wideness \newcount\leftwidth \newcount\rightwidth
\newcount\highness \newcount\topheight \newcount\botheight

\def\ratchet#1#2{\ifnum#1<#2 \global #1=#2 \fi}

\def\putbox(#1,#2)#3{%
\horsize{\wideness}{#3} \divide\wideness by 2
{\advance\wideness by #1 \ratchet{\rightwidth}{\wideness}}
{\advance\wideness by -#1 \ratchet{\leftwidth}{\wideness}}
\vertsize{\highness}{#3} \divide\highness by 2
{\advance\highness by #2 \ratchet{\topheight}{\highness}}
{\advance\highness by -#2 \ratchet{\botheight}{\highness}}
\put(#1,#2){\makebox(0,0){$#3$}}}

\def\putlbox(#1,#2)#3{%
\horsize{\wideness}{#3}
{\advance\wideness by #1 \ratchet{\rightwidth}{\wideness}}
{\ratchet{\leftwidth}{-#1}}
\vertsize{\highness}{#3} \divide\highness by 2
{\advance\highness by #2 \ratchet{\topheight}{\highness}}
{\advance\highness by -#2 \ratchet{\botheight}{\highness}}
\put(#1,#2){\makebox(0,0)[l]{$#3$}}}

\def\putrbox(#1,#2)#3{%
\horsize{\wideness}{#3}
{\ratchet{\rightwidth}{#1}}
{\advance\wideness by -#1 \ratchet{\leftwidth}{\wideness}}
\vertsize{\highness}{#3} \divide\highness by 2
{\advance\highness by #2 \ratchet{\topheight}{\highness}}
{\advance\highness by -#2 \ratchet{\botheight}{\highness}}
\put(#1,#2){\makebox(0,0)[r]{$#3$}}}

\def\adjust[#1]{} 

\newcount \coefa
\newcount \coefb
\newcount \coefc
\newcount\tempcounta
\newcount\tempcountb
\newcount\tempcountc
\newcount\tempcountd
\newcount\xext
\newcount\yext
\newcount\xoff
\newcount\yoff
\newcount\gap%
\newcount\arrowtypea
\newcount\arrowtypeb
\newcount\arrowtypec
\newcount\arrowtyped
\newcount\arrowtypee
\newcount\height
\newcount\width
\newcount\xpos
\newcount\ypos
\newcount\run
\newcount\rise
\newcount\arrowlength
\newcount\halflength
\newcount\arrowtype
\newdimen\tempdimen
\newdimen\xlen
\newdimen\ylen
\newsavebox{\tempboxa}%
\newsavebox{\tempboxb}%
\newsavebox{\tempboxc}%

\newdimen\w@dth

\def\setw@dth#1#2{\setbox\z@\hbox{\m@th$#1$}\w@dth=\wd\z@
\setbox\@ne\hbox{\m@th$#2$}\ifnum\w@dth<\wd\@ne \w@dth=\wd\@ne \fi
\advance\w@dth by 1.2em}

\def\t@^#1_#2{\allowbreak\def\n@one{#1}\def\n@two{#2}\mathrel
{\setw@dth{#1}{#2}
\mathop{\hbox to \w@dth{\rightarrowfill}}\limits
\ifx\n@one\empty\else ^{\box\z@}\fi
\ifx\n@two\empty\else _{\box\@ne}\fi}}
\def\t@@^#1{\@ifnextchar_{\t@^{#1}}{\t@^{#1}_{}}}
\def\to{\@ifnextchar^{\t@@}{\t@@^{}}}

\def\t@left^#1_#2{\def\n@one{#1}\def\n@two{#2}\mathrel{\setw@dth{#1}{#2}
\mathop{\hbox to \w@dth{\leftarrowfill}}\limits
\ifx\n@one\empty\else ^{\box\z@}\fi
\ifx\n@two\empty\else _{\box\@ne}\fi}}
\def\t@@left^#1{\@ifnextchar_{\t@left^{#1}}{\t@left^{#1}_{}}}
\def\toleft{\@ifnextchar^{\t@@left}{\t@@left^{}}}

\def\two@^#1_#2{\allowbreak
\def\n@one{#1}\def\n@two{#2}\mathrel{\setw@dth{#1}{#2}
\mathop{\vcenter{\lineskip\z@\baselineskip\z@
                 \hbox to \w@dth{\rightarrowfill}%
                 \hbox to \w@dth{\rightarrowfill}}%
       }\limits
\ifx\n@one\empty\else ^{\box\z@}\fi
\ifx\n@two\empty\else _{\box\@ne}\fi}}
\def\tw@@^#1{\@ifnextchar _{\two@^{#1}}{\two@^{#1}_{}}}
\def\two{\@ifnextchar ^{\tw@@}{\tw@@^{}}}

\def\tofr@^#1_#2{\def\n@one{#1}\def\n@two{#2}\mathrel{\setw@dth{#1}{#2}
\mathop{\vcenter{\hbox to \w@dth{\rightarrowfill}\kern-1.7ex
                 \hbox to \w@dth{\leftarrowfill}}%
       }\limits
\ifx\n@one\empty\else ^{\box\z@}\fi
\ifx\n@two\empty\else _{\box\@ne}\fi}}
\def\t@fr@^#1{\@ifnextchar_ {\tofr@^{#1}}{\tofr@^{#1}_{}}}
\def\tofro{\@ifnextchar^ {\t@fr@}{\t@fr@^{}}}

\def\mon{\mathop{\m@th\hbox to
      14.6\P@{\lasyb\char'51\hskip-2.1\P@$\arrext$\hss
$\mathord\rightarrow$}}\limits} 
\def\leftmono{\mathrel{\m@th\hbox to
14.6\P@{$\mathord\leftarrow$\hss$\arrext$\hskip-2.1\P@\lasyb\char'50%
}}\limits} 
\mathchardef\arrext="0200       

\def\settypes(#1,#2,#3){\arrowtypea#1 \arrowtypeb#2 \arrowtypec#3}
\def\settoheight#1#2{\setbox\@tempboxa\hbox{#2}#1\ht\@tempboxa\relax}%
\def\settodepth#1#2{\setbox\@tempboxa\hbox{#2}#1\dp\@tempboxa\relax}%
\def\settokens`#1`#2`#3`#4`{%
     \def\tokena{#1}\def\tokenb{#2}\def\tokenc{#3}\def\tokend{#4}}
\def\setsqparms[#1`#2`#3`#4;#5`#6]{%
\arrowtypea #1
\arrowtypeb #2
\arrowtypec #3
\arrowtyped #4
\width #5
\height #6
}
\def\setpos(#1,#2){\xpos=#1 \ypos#2}

\def\settriparms[#1`#2`#3;#4]{\settripairparms[#1`#2`#3`1`1;#4]}%

\def\settripairparms[#1`#2`#3`#4`#5;#6]{%
\arrowtypea #1
\arrowtypeb #2
\arrowtypec #3
\arrowtyped #4
\arrowtypee #5
\width #6
\height #6
}

\def\resetparms{\settripairparms[1`1`1`1`1;500]\width 500}

\resetparms

\def\mvector(#1,#2)#3{
\put(0,0){\vector(#1,#2){#3}}%
\put(0,0){\vector(#1,#2){26}}%
}
\def\evector(#1,#2)#3{{
\arrowlength #3
\put(0,0){\vector(#1,#2){\arrowlength}}%
\advance \arrowlength by-30
\put(0,0){\vector(#1,#2){\arrowlength}}%
}}

\def\horsize#1#2{%
\settowidth{\tempdimen}{$#2$}%
#1=\tempdimen
\divide #1 by\unitlength
}

\def\vertsize#1#2{%
\settoheight{\tempdimen}{$#2$}%
#1=\tempdimen
\settodepth{\tempdimen}{$#2$}%
\advance #1 by\tempdimen
\divide #1 by\unitlength
}

\def\putvector(#1,#2)(#3,#4)#5#6{{%
\ifnum3<\arrowtype
\putdashvector(#1,#2)(#3,#4)#5\arrowtype
\else
\ifnum\arrowtype<-3
\putdashvector(#1,#2)(#3,#4)#5\arrowtype
\else
\xpos=#1
\ypos=#2
\run=#3
\rise=#4
\arrowlength=#5
\ifnum \arrowtype<0
    \ifnum \run=0
        \advance \ypos by-\arrowlength
    \else
        \tempcounta \arrowlength
        \multiply \tempcounta by\rise
        \divide \tempcounta by\run
        \ifnum\run>0
            \advance \xpos by\arrowlength
            \advance \ypos by\tempcounta
        \else
            \advance \xpos by-\arrowlength
            \advance \ypos by-\tempcounta
        \fi
    \fi
    \multiply \arrowtype by-1
    \multiply \rise by-1
    \multiply \run by-1
\fi
\ifcase \arrowtype
\or \put(\xpos,\ypos){\vector(\run,\rise){\arrowlength}}%
\or \put(\xpos,\ypos){\mvector(\run,\rise)\arrowlength}%
\or \put(\xpos,\ypos){\evector(\run,\rise){\arrowlength}}%
\fi\fi\fi
}}

\def\putsplitvector(#1,#2)#3#4{
\xpos #1
\ypos #2
\arrowtype #4
\halflength #3
\arrowlength #3
\gap 140
\advance \halflength by-\gap
\divide \halflength by2
\ifnum\arrowtype>0
   \ifcase \arrowtype
   \or \put(\xpos,\ypos){\line(0,-1){\halflength}}%
       \advance\ypos by-\halflength
       \advance\ypos by-\gap
       \put(\xpos,\ypos){\vector(0,-1){\halflength}}%
   \or \put(\xpos,\ypos){\line(0,-1)\halflength}%
       \put(\xpos,\ypos){\vector(0,-1)3}%
       \advance\ypos by-\halflength
       \advance\ypos by-\gap
       \put(\xpos,\ypos){\vector(0,-1){\halflength}}%
   \or \put(\xpos,\ypos){\line(0,-1)\halflength}%
       \advance\ypos by-\halflength
       \advance\ypos by-\gap
       \put(\xpos,\ypos){\evector(0,-1){\halflength}}%
   \fi
\else \arrowtype=-\arrowtype
   \ifcase\arrowtype
   \or \advance \ypos by-\arrowlength
       \put(\xpos,\ypos){\line(0,1){\halflength}}%
       \advance\ypos by\halflength
       \advance\ypos by\gap
       \put(\xpos,\ypos){\vector(0,1){\halflength}}%
   \or \advance \ypos by-\arrowlength
       \put(\xpos,\ypos){\line(0,1)\halflength}%
       \put(\xpos,\ypos){\vector(0,1)3}%
       \advance\ypos by\halflength
       \advance\ypos by\gap
       \put(\xpos,\ypos){\vector(0,1){\halflength}}%
   \or \advance \ypos by-\arrowlength
       \put(\xpos,\ypos){\line(0,1)\halflength}%
       \advance\ypos by\halflength
       \advance\ypos by\gap
       \put(\xpos,\ypos){\evector(0,1){\halflength}}%
   \fi
\fi
}

\def\putmorphism(#1)(#2,#3)[#4`#5`#6]#7#8#9{{%
\run #2
\rise #3
\ifnum\rise=0
  \puthmorphism(#1)[#4`#5`#6]{#7}{#8}#9%
\else\ifnum\run=0
  \putvmorphism(#1)[#4`#5`#6]{#7}{#8}#9%
\else
\setpos(#1)%
\arrowlength #7
\arrowtype #8
\ifnum\run=0
\else\ifnum\rise=0
\else
\ifnum\run>0
    \coefa=1
\else
   \coefa=-1
\fi
\ifnum\arrowtype>0
   \coefb=0
   \coefc=-1
\else
   \coefb=\coefa
   \coefc=1
   \arrowtype=-\arrowtype
\fi
\width=2
\multiply \width by\run
\divide \width by\rise
\ifnum \width<0  \width=-\width\fi
\advance\width by60
\if l#9 \width=-\width\fi
\putbox(\xpos,\ypos){#4}
{\multiply \coefa by\arrowlength
\advance\xpos by\coefa
\multiply \coefa by\rise
\divide \coefa by\run
\advance \ypos by\coefa
\putbox(\xpos,\ypos){#5} }%
{\multiply \coefa by\arrowlength
\divide \coefa by2
\advance \xpos by\coefa
\advance \xpos by\width
\multiply \coefa by\rise
\divide \coefa by\run
\advance \ypos by\coefa
\if l#9%
   \putrbox(\xpos,\ypos){#6}%
\else\if r#9%
   \putlbox(\xpos,\ypos){#6}%
\fi\fi }%
{\multiply \rise by-\coefc
\multiply \run by-\coefc
\multiply \coefb by\arrowlength
\advance \xpos by\coefb
\multiply \coefb by\rise
\divide \coefb by\run
\advance \ypos by\coefb
\multiply \coefc by70
\advance \ypos by\coefc
\multiply \coefc by\run
\divide \coefc by\rise
\advance \xpos by\coefc
\multiply \coefa by140
\multiply \coefa by\run
\divide \coefa by\rise
\advance \arrowlength by\coefa
\ifcase\arrowtype
\or \put(\xpos,\ypos){\vector(\run,\rise){\arrowlength}}%
\or \put(\xpos,\ypos){\mvector(\run,\rise){\arrowlength}}%
\or \put(\xpos,\ypos){\evector(\run,\rise){\arrowlength}}%
\fi}\fi\fi\fi\fi}}

\newcount\numbdashes \newcount\lengthdash \newcount\increment

\def\howmanydashes{
\numbdashes=\arrowlength \lengthdash=40
\divide\numbdashes by \lengthdash
\lengthdash=\arrowlength
\divide\lengthdash by \numbdashes
\increment=\lengthdash
\multiply\lengthdash by 3
\divide\lengthdash by 5
}

\def\putdashvector(#1)(#2,#3)#4#5{%
\ifnum#3=0 \putdashhvector(#1){#4}#5
\else
\ifnum#2=0
\putdashvvector(#1){#4}#5\fi\fi}

\def\putdashhvector(#1,#2)#3#4{{%
\arrowlength=#3 \howmanydashes
\multiput(#1,#2)(\increment,0){\numbdashes}%
{\vrule height .4pt width \lengthdash\unitlength}
\arrowtype=#4 \xpos=#1
\ifnum\arrowtype<0 \advance\arrowtype by 7 \fi
\ifcase\arrowtype
\or \advance\xpos by 10
    \put(\xpos,#2){\vector(-1,0){\lengthdash}}
    \advance\xpos by 40
    \put(\xpos,#2){\vector(-1,0){\lengthdash}}
\or \advance \xpos by 10
    \put(\xpos,#2){\vector(-1,0){\lengthdash}}
    \advance\xpos by  \arrowlength
    \advance\xpos by  -50
    \put(\xpos,#2){\vector(-1,0){\lengthdash}}
\or \advance\xpos by 10
    \put(\xpos,#2){\vector(-1,0){\lengthdash}}
\or \advance\xpos by \arrowlength
    \advance\xpos by -\lengthdash
    \put(\xpos,#2){\vector(1,0){\lengthdash}}
\or {\advance\xpos by 10
    \put(\xpos,#2){\vector(1,0){\lengthdash}}}
    \advance\xpos by \arrowlength
    \advance\xpos by -\lengthdash
    \put(\xpos,#2){\vector(1,0){\lengthdash}}
\or \advance\xpos by \arrowlength
    \advance\xpos by -\lengthdash
    \put(\xpos,#2){\vector(1,0){\lengthdash}}
    \advance\xpos by -40
    \put(\xpos,#2){\vector(1,0){\lengthdash}}
   \fi
}}

\def\putdashvvector(#1,#2)#3#4{{%
\arrowlength=#3 \howmanydashes
\ypos=#2 \advance\ypos by -\arrowlength
\multiput(#1,#2)(0,\increment){\numbdashes}%
    {\vrule width .4pt height \lengthdash\unitlength}
\arrowtype=#4 \ypos=#2
\ifnum\arrowtype<0 \advance\arrowtype by 7 \fi
\ifcase\arrowtype
\or \advance\ypos by \arrowlength \advance\ypos by -40
    \put(#1,\ypos){\vector(0,1){\lengthdash}}
    \advance\ypos by -40
    \put(#1,\ypos){\vector(0,1){\lengthdash}}
\or \advance\ypos by 10
    \put(#1,\ypos){\vector(0,1){\lengthdash}}
    \advance\ypos by \arrowlength \advance\ypos by -40
    \put(#1,\ypos){\vector(0,1){\lengthdash}}
\or \advance\ypos by \arrowlength \advance\ypos by -40
    \put(#1,\ypos){\vector(0,1){\lengthdash}}
\or \advance\ypos by 10
    \put(#1,\ypos){\vector(0,-1){\lengthdash}}
\or \advance\ypos by 10
    \put(#1,\ypos){\vector(0,-1){\lengthdash}}
    \advance\ypos by \arrowlength \advance\ypos by -40
    \put(#1,\ypos){\vector(0,-1){\lengthdash}}
\or \advance\ypos by 10
    \put(#1,\ypos){\vector(0,-1){\lengthdash}}
    \advance\ypos by 40
    \put(#1,\ypos){\vector(0,-1){\lengthdash}}
\fi
}}

\def\puthmorphism(#1,#2)[#3`#4`#5]#6#7#8{{%
\xpos #1
\ypos #2
\width #6
\arrowlength #6
\arrowtype=#7
\putbox(\xpos,\ypos){#3\vphantom{#4}}%
{\advance \xpos by\arrowlength
\putbox(\xpos,\ypos){\vphantom{#3}#4}}%
\horsize{\tempcounta}{#3}%
\horsize{\tempcountb}{#4}%
\divide \tempcounta by2
\divide \tempcountb by2
\advance \tempcounta by30
\advance \tempcountb by30
\advance \xpos by\tempcounta
\advance \arrowlength by-\tempcounta
\advance \arrowlength by-\tempcountb
\putvector(\xpos,\ypos)(1,0)\arrowlength\arrowtype
\divide \arrowlength by2
\advance \xpos by\arrowlength
\vertsize{\tempcounta}{#5}%
\divide\tempcounta by2
\advance \tempcounta by20
\if a#8 %
   \advance \ypos by\tempcounta
   \putbox(\xpos,\ypos){#5}%
\else
   \advance \ypos by-\tempcounta
   \putbox(\xpos,\ypos){#5}%
\fi}}

\def\putvmorphism(#1,#2)[#3`#4`#5]#6#7#8{{%
\xpos #1
\ypos #2
\arrowlength #6
\arrowtype #7
\settowidth{\xlen}{$#5$}%
\putbox(\xpos,\ypos){#3}%
{\advance \ypos by-\arrowlength
\putbox(\xpos,\ypos){#4}}%
{\advance\arrowlength by-140
\advance \ypos by-70
\ifdim\xlen>0pt
   \if m#8%
      \putsplitvector(\xpos,\ypos)\arrowlength\arrowtype
   \else
   \putvector(\xpos,\ypos)(0,-1)\arrowlength\arrowtype
   \fi
\else
   \putvector(\xpos,\ypos)(0,-1)\arrowlength\arrowtype
\fi}%
\ifdim\xlen>0pt
   \divide \arrowlength by2
   \advance\ypos by-\arrowlength
   \if l#8%
      \advance \xpos by-40
      \putrbox(\xpos,\ypos){#5}%
   \else\if r#8%
      \advance \xpos by40
      \putlbox(\xpos,\ypos){#5}%
   \else
      \putbox(\xpos,\ypos){#5}%
   \fi\fi
\fi
}}

\def\putsquarep<#1>(#2)[#3;#4`#5`#6`#7]{{%
\setsqparms[#1]%
\setpos(#2)%
\settokens`#3`%
\puthmorphism(\xpos,\ypos)[\tokenc`\tokend`{#7}]{\width}{\arrowtyped}b%
\advance\ypos by \height
\puthmorphism(\xpos,\ypos)[\tokena`\tokenb`{#4}]{\width}{\arrowtypea}a%
\putvmorphism(\xpos,\ypos)[``{#5}]{\height}{\arrowtypeb}l%
\advance\xpos by \width
\putvmorphism(\xpos,\ypos)[``{#6}]{\height}{\arrowtypec}r%
}}

\def\putsquare{\@ifnextchar <{\putsquarep}{\putsquarep%
   <\arrowtypea`\arrowtypeb`\arrowtypec`\arrowtyped;\width`\height>}}
\def\square{\@ifnextchar< {\squarep}{\squarep
   <\arrowtypea`\arrowtypeb`\arrowtypec`\arrowtyped;\width`\height>}}
\def\squarep<#1>[#2`#3`#4`#5;#6`#7`#8`#9]{{
\setsqparms[#1]
\diagram
\putsquarep<\arrowtypea`\arrowtypeb`\arrowtypec`
\arrowtyped;\width`\height>
(0,0)[#2`#3`#4`{#5};#6`#7`#8`{#9}]
\enddiagram
}}                                                 
\def\putptrianglep<#1>(#2,#3)[#4`#5`#6;#7`#8`#9]{{%
\settriparms[#1]%
\xpos=#2 \ypos=#3
\advance\ypos by \height
\puthmorphism(\xpos,\ypos)[#4`#5`{#7}]{\height}{\arrowtypea}a%
\putvmorphism(\xpos,\ypos)[`#6`{#8}]{\height}{\arrowtypeb}l%
\advance\xpos by\height
\putmorphism(\xpos,\ypos)(-1,-1)[``{#9}]{\height}{\arrowtypec}r%
}}

\def\putptriangle{\@ifnextchar <{\putptrianglep}{\putptrianglep
   <\arrowtypea`\arrowtypeb`\arrowtypec;\height>}}
\def\ptriangle{\@ifnextchar <{\ptrianglep}{\ptrianglep
   <\arrowtypea`\arrowtypeb`\arrowtypec;\height>}}
\def\ptrianglep<#1>[#2`#3`#4;#5`#6`#7]{{
\settriparms[#1]
\diagram
\putptrianglep<\arrowtypea`\arrowtypeb`
\arrowtypec;\height>
(0,0)[#2`#3`#4;#5`#6`{#7}]
\enddiagram
}}                                            

\def\putqtrianglep<#1>(#2,#3)[#4`#5`#6;#7`#8`#9]{{%
\settriparms[#1]%
\xpos=#2 \ypos=#3
\advance\ypos by\height
\puthmorphism(\xpos,\ypos)[#4`#5`{#7}]{\height}{\arrowtypea}a%
\putmorphism(\xpos,\ypos)(1,-1)[``{#8}]{\height}{\arrowtypeb}l%
\advance\xpos by\height
\putvmorphism(\xpos,\ypos)[`#6`{#9}]{\height}{\arrowtypec}r%
}}

\def\putqtriangle{\@ifnextchar <{\putqtrianglep}{\putqtrianglep
   <\arrowtypea`\arrowtypeb`\arrowtypec;\height>}}
\def\qtriangle{\@ifnextchar <{\qtrianglep}{\qtrianglep
   <\arrowtypea`\arrowtypeb`\arrowtypec;\height>}}
\def\qtrianglep<#1>[#2`#3`#4;#5`#6`#7]{{
\settriparms[#1]
\width=\height                                
\diagram
\putqtrianglep<\arrowtypea`\arrowtypeb`
\arrowtypec;\height>
(0,0)[#2`#3`#4;#5`#6`{#7}]
\enddiagram
}}

\def\putdtrianglep<#1>(#2,#3)[#4`#5`#6;#7`#8`#9]{{%
\settriparms[#1]%
\xpos=#2 \ypos=#3
\puthmorphism(\xpos,\ypos)[#5`#6`{#9}]{\height}{\arrowtypec}b%
\advance\xpos by \height \advance\ypos by\height
\putmorphism(\xpos,\ypos)(-1,-1)[``{#7}]{\height}{\arrowtypea}l%
\putvmorphism(\xpos,\ypos)[#4``{#8}]{\height}{\arrowtypeb}r%
}}

\def\putdtriangle{\@ifnextchar <{\putdtrianglep}{\putdtrianglep
   <\arrowtypea`\arrowtypeb`\arrowtypec;\height>}}
\def\dtriangle{\@ifnextchar <{\dtrianglep}{\dtrianglep
   <\arrowtypea`\arrowtypeb`\arrowtypec;\height>}}
\def\dtrianglep<#1>[#2`#3`#4;#5`#6`#7]{{
\settriparms[#1]
\width=\height                                
\diagram
\putdtrianglep<\arrowtypea`\arrowtypeb`
\arrowtypec;\height>
(0,0)[#2`#3`#4;#5`#6`{#7}]
\enddiagram
}}

\def\putbtrianglep<#1>(#2,#3)[#4`#5`#6;#7`#8`#9]{{%
\settriparms[#1]%
\xpos=#2 \ypos=#3
\puthmorphism(\xpos,\ypos)[#5`#6`{#9}]{\height}{\arrowtypec}b%
\advance\ypos by\height
\putmorphism(\xpos,\ypos)(1,-1)[``{#8}]{\height}{\arrowtypeb}r%
\putvmorphism(\xpos,\ypos)[#4``{#7}]{\height}{\arrowtypea}l%
}}

\def\putbtriangle{\@ifnextchar <{\putbtrianglep}{\putbtrianglep
   <\arrowtypea`\arrowtypeb`\arrowtypec;\height>}}
\def\btriangle{\@ifnextchar <{\btrianglep}{\btrianglep
   <\arrowtypea`\arrowtypeb`\arrowtypec;\height>}}
\def\btrianglep<#1>[#2`#3`#4;#5`#6`#7]{{
\settriparms[#1]
\width=\height                               
\diagram
\putbtrianglep<\arrowtypea`\arrowtypeb`
\arrowtypec;\height>
(0,0)[#2`#3`#4;#5`#6`{#7}]
\enddiagram
}}

\def\putAtrianglep<#1>(#2,#3)[#4`#5`#6;#7`#8`#9]{{%
\settriparms[#1]%
\xpos=#2 \ypos=#3
{\multiply \height by2
\puthmorphism(\xpos,\ypos)[#5`#6`{#9}]{\height}{\arrowtypec}b}%
\advance\xpos by\height \advance\ypos by\height
\putmorphism(\xpos,\ypos)(-1,-1)[#4``{#7}]{\height}{\arrowtypea}l%
\putmorphism(\xpos,\ypos)(1,-1)[``{#8}]{\height}{\arrowtypeb}r%
}}

\def\putAtriangle{\@ifnextchar <{\putAtrianglep}{\putAtrianglep
   <\arrowtypea`\arrowtypeb`\arrowtypec;\height>}}
\def\Atriangle{\@ifnextchar <{\Atrianglep}{\Atrianglep
   <\arrowtypea`\arrowtypeb`\arrowtypec;\height>}}
\def\Atrianglep<#1>[#2`#3`#4;#5`#6`#7]{{
\settriparms[#1]
\width=\height                                     
\diagram
\putAtrianglep<\arrowtypea`\arrowtypeb`
\arrowtypec;\height>
(0,0)[#2`#3`#4;#5`#6`{#7}]
\enddiagram
}}

\def\putAtrianglepairp<#1>(#2)[#3;#4`#5`#6`#7`#8]{{%
\settripairparms[#1]%
\setpos(#2)%
\settokens`#3`%
\puthmorphism(\xpos,\ypos)[\tokenb`\tokenc`{#7}]{\height}{\arrowtyped}b%
\advance\xpos by\height
\puthmorphism(\xpos,\ypos)[\phantom{\tokenc}`\tokend`{#8}]%
{\height}{\arrowtypee}b%
\advance\ypos by\height
\putmorphism(\xpos,\ypos)(-1,-1)[\tokena``{#4}]{\height}{\arrowtypea}l%
\putvmorphism(\xpos,\ypos)[``{#5}]{\height}{\arrowtypeb}m%
\putmorphism(\xpos,\ypos)(1,-1)[``{#6}]{\height}{\arrowtypec}r%
}}

\def\putAtrianglepair{\@ifnextchar <{\putAtrianglepairp}{\putAtrianglepairp%
   <\arrowtypea`\arrowtypeb`\arrowtypec`\arrowtyped`\arrowtypee;\height>}}
\def\Atrianglepair{\@ifnextchar <{\Atrianglepairp}{\Atrianglepairp%
   <\arrowtypea`\arrowtypeb`\arrowtypec`\arrowtyped`\arrowtypee;\height>}}

\def\Atrianglepairp<#1>[#2;#3`#4`#5`#6`#7]{{
\settripairparms[#1]
\settokens`#2`
\width=\height                                
\diagram
\putAtrianglepairp                            
<\arrowtypea`\arrowtypeb`\arrowtypec`
\arrowtyped`\arrowtypee;\height>
(0,0)[{#2};#3`#4`#5`#6`{#7}]
\enddiagram
}}

\def\putVtrianglep<#1>(#2,#3)[#4`#5`#6;#7`#8`#9]{{%
\settriparms[#1]%
\xpos=#2 \ypos=#3
\advance\ypos by\height
{\multiply\height by2
\puthmorphism(\xpos,\ypos)[#4`#5`{#7}]{\height}{\arrowtypea}a}%
\putmorphism(\xpos,\ypos)(1,-1)[`#6`{#8}]{\height}{\arrowtypeb}l%
\advance\xpos by\height
\advance\xpos by\height
\putmorphism(\xpos,\ypos)(-1,-1)[``{#9}]{\height}{\arrowtypec}r%
}}

\def\putVtriangle{\@ifnextchar <{\putVtrianglep}{\putVtrianglep
   <\arrowtypea`\arrowtypeb`\arrowtypec;\height>}}
\def\Vtriangle{\@ifnextchar <{\Vtrianglep}{\Vtrianglep
   <\arrowtypea`\arrowtypeb`\arrowtypec;\height>}}
\def\Vtrianglep<#1>[#2`#3`#4;#5`#6`#7]{{
\settriparms[#1]
\width=\height                                 
\diagram
\putVtrianglep<\arrowtypea`\arrowtypeb`
\arrowtypec;\height>
(0,0)[#2`#3`#4;#5`#6`{#7}]
\enddiagram
}}

\def\putVtrianglepairp<#1>(#2)[#3;#4`#5`#6`#7`#8]{{
\settripairparms[#1]%
\setpos(#2)%
\settokens`#3`%
\advance\ypos by\height
\putmorphism(\xpos,\ypos)(1,-1)[`\tokend`{#6}]{\height}{\arrowtypec}l%
\puthmorphism(\xpos,\ypos)[\tokena`\tokenb`{#4}]{\height}{\arrowtypea}a%
\advance\xpos by\height
\puthmorphism(\xpos,\ypos)[\phantom{\tokenb}`\tokenc`{#5}]%
{\height}{\arrowtypeb}a%
\putvmorphism(\xpos,\ypos)[``{#7}]{\height}{\arrowtyped}m%
\advance\xpos by\height
\putmorphism(\xpos,\ypos)(-1,-1)[``{#8}]{\height}{\arrowtypee}r%
}}

\def\putVtrianglepair{\@ifnextchar <{\putVtrianglepairp}{\putVtrianglepairp%
    <\arrowtypea`\arrowtypeb`\arrowtypec`\arrowtyped`\arrowtypee;\height>}}
\def\Vtrianglepair{\@ifnextchar <{\Vtrianglepairp}{\Vtrianglepairp%
    <\arrowtypea`\arrowtypeb`\arrowtypec`\arrowtyped`\arrowtypee;\height>}}
\def\Vtrianglepairp<#1>[#2;#3`#4`#5`#6`#7]{{
\settripairparms[#1]
\settokens`#2`
\diagram
\putVtrianglepairp                             
<\arrowtypea`\arrowtypeb`\arrowtypec`
\arrowtyped`\arrowtypee;\height>
(0,0)[{#2};#3`#4`#5`#6`{#7}]
\enddiagram
}}

\def\putCtrianglep<#1>(#2,#3)[#4`#5`#6;#7`#8`#9]{{%
\settriparms[#1]%
\xpos=#2 \ypos=#3
\advance\ypos by\height
\putmorphism(\xpos,\ypos)(1,-1)[``{#9}]{\height}{\arrowtypec}l%
\advance\xpos by\height
\advance\ypos by\height
\putmorphism(\xpos,\ypos)(-1,-1)[#4`#5`{#7}]{\height}{\arrowtypea}l%
{\multiply\height by 2
\putvmorphism(\xpos,\ypos)[`#6`{#8}]{\height}{\arrowtypeb}r}%
}}

\def\putCtriangle{\@ifnextchar <{\putCtrianglep}{\putCtrianglep
    <\arrowtypea`\arrowtypeb`\arrowtypec;\height>}}
\def\Ctriangle{\@ifnextchar <{\Ctrianglep}{\Ctrianglep
    <\arrowtypea`\arrowtypeb`\arrowtypec;\height>}}
\def\Ctrianglep<#1>[#2`#3`#4;#5`#6`#7]{{
\settriparms[#1]
\width=\height                               
\diagram
\putCtrianglep<\arrowtypea`\arrowtypeb`
\arrowtypec;\height>
(0,0)[#2`#3`#4;#5`#6`{#7}]
\enddiagram
}}                                           
\def\putDtrianglep<#1>(#2,#3)[#4`#5`#6;#7`#8`#9]{{%
\settriparms[#1]%
\xpos=#2 \ypos=#3
\advance\xpos by\height \advance\ypos by\height
\putmorphism(\xpos,\ypos)(-1,-1)[``{#9}]{\height}{\arrowtypec}r%
\advance\xpos by-\height \advance\ypos by\height
\putmorphism(\xpos,\ypos)(1,-1)[`#5`{#8}]{\height}{\arrowtypeb}r%
{\multiply\height by 2
\putvmorphism(\xpos,\ypos)[#4`#6`{#7}]{\height}{\arrowtypea}l}%
}}

\def\putDtriangle{\@ifnextchar <{\putDtrianglep}{\putDtrianglep
    <\arrowtypea`\arrowtypeb`\arrowtypec;\height>}}
\def\Dtriangle{\@ifnextchar <{\Dtrianglep}{\Dtrianglep
   <\arrowtypea`\arrowtypeb`\arrowtypec;\height>}}
\def\Dtrianglep<#1>[#2`#3`#4;#5`#6`#7]{{
\settriparms[#1]
\width=\height                              
\diagram
\putDtrianglep<\arrowtypea`\arrowtypeb`
\arrowtypec;\height>
(0,0)[#2`#3`#4;#5`#6`{#7}]
\enddiagram
}}                                          
\def\setrecparms[#1`#2]{\width=#1 \height=#2}%

\def\recursep<#1`#2>[#3;#4`#5`#6`#7`#8]{{\m@th
\width=#1 \height=#2
\settokens`#3`
\settowidth{\tempdimen}{$\tokena$}
\ifdim\tempdimen=0pt
  \savebox{\tempboxa}{\hbox{$\tokenb$}}%
  \savebox{\tempboxb}{\hbox{$\tokend$}}%
  \savebox{\tempboxc}{\hbox{$#6$}}%
\else
  \savebox{\tempboxa}{\hbox{$\hbox{$\tokena$}\times\hbox{$\tokenb$}$}}%
  \savebox{\tempboxb}{\hbox{$\hbox{$\tokena$}\times\hbox{$\tokend$}$}}%
  \savebox{\tempboxc}{\hbox{$\hbox{$\tokena$}\times\hbox{$#6$}$}}%
\fi
\ypos=\height
\divide\ypos by 2
\xpos=\ypos
\advance\xpos by \width
\bfig
\putCtrianglep<-1`1`1;\ypos>(0,0)[`\tokenc`;#5`#6`{#7}]%
\puthmorphism(\ypos,0)[\tokend`\usebox{\tempboxb}`{#8}]{\width}{-1}b%
\puthmorphism(\ypos,\height)[\tokenb`\usebox{\tempboxa}`{#4}]{\width}{-1}a%
\advance\ypos by \width
\putvmorphism(\ypos,\height)[``\usebox{\tempboxc}]{\height}1r%
\efig
}}

\def\recurse{\@ifnextchar <{\recursep}{\recursep<\width`\height>}}

\def\puttwohmorphisms(#1,#2)[#3`#4;#5`#6]#7#8#9{{%
\puthmorphism(#1,#2)[#3`#4`]{#7}0a
\ypos=#2
\advance\ypos by 20
\puthmorphism(#1,\ypos)[\phantom{#3}`\phantom{#4}`#5]{#7}{#8}a
\advance\ypos by -40
\puthmorphism(#1,\ypos)[\phantom{#3}`\phantom{#4}`#6]{#7}{#9}b
}}

\def\puttwovmorphisms(#1,#2)[#3`#4;#5`#6]#7#8#9{{%
\putvmorphism(#1,#2)[#3`#4`]{#7}0a
\xpos=#1
\advance\xpos by -20
\putvmorphism(\xpos,#2)[\phantom{#3}`\phantom{#4}`#5]{#7}{#8}l
\advance\xpos by 40
\putvmorphism(\xpos,#2)[\phantom{#3}`\phantom{#4}`#6]{#7}{#9}r
}}

\def\puthcoequalizer(#1)[#2`#3`#4;#5`#6`#7]#8#9{{%
\setpos(#1)%
\puttwohmorphisms(\xpos,\ypos)[#2`#3;#5`#6]{#8}11%
\advance\xpos by #8
\puthmorphism(\xpos,\ypos)[\phantom{#3}`#4`#7]{#8}1{#9}
}}

\def\putvcoequalizer(#1)[#2`#3`#4;#5`#6`#7]#8#9{{%
\setpos(#1)%
\puttwovmorphisms(\xpos,\ypos)[#2`#3;#5`#6]{#8}11%
\advance\ypos by -#8
\putvmorphism(\xpos,\ypos)[\phantom{#3}`#4`#7]{#8}1{#9}
}}

\def\putthreehmorphisms(#1)[#2`#3;#4`#5`#6]#7(#8)#9{{%
\setpos(#1) \settypes(#8)
\if a#9 %
     \vertsize{\tempcounta}{#5}%
     \vertsize{\tempcountb}{#6}%
     \ifnum \tempcounta<\tempcountb \tempcounta=\tempcountb \fi
\else
     \vertsize{\tempcounta}{#4}%
     \vertsize{\tempcountb}{#5}%
     \ifnum \tempcounta<\tempcountb \tempcounta=\tempcountb \fi
\fi
\advance \tempcounta by 60
\puthmorphism(\xpos,\ypos)[#2`#3`#5]{#7}{\arrowtypeb}{#9}
\advance\ypos by \tempcounta
\puthmorphism(\xpos,\ypos)[\phantom{#2}`\phantom{#3}`#4]{#7}{\arrowtypea}{#9}
\advance\ypos by -\tempcounta \advance\ypos by -\tempcounta
\puthmorphism(\xpos,\ypos)[\phantom{#2}`\phantom{#3}`#6]{#7}{\arrowtypec}{#9}
}}

\def\setarrowtoks[#1`#2`#3`#4`#5`#6]{%
\def\toka{#1}
\def\tokb{#2}
\def\tokc{#3}
\def\tokd{#4}
\def\toke{#5}
\def\tokf{#6}
}
\def\hex{\@ifnextchar <{\hexp}{\hexp<1000`400>}}
\def\hexp<#1`#2>[#3`#4`#5`#6`#7`#8;#9]{%
\setarrowtoks[#9]
\yext=#2 \advance \yext by #2
\xext=#1 \advance\xext by \yext
\bfig
\putCtriangle<-1`0`1;#2>(0,0)[`#5`;\tokb``\tokd]
\xext=#1 \yext=#2 \advance \yext by #2
\putsquare<1`0`0`1;\xext`\yext>(#2,0)[#3`#4`#7`#8;\toka```\tokf]
\advance \xext by #2
\putDtriangle<0`1`-1;#2>(\xext,0)[`#6`;`\tokc`\toke]
\efig
}
\makeatother

\begin{document}
\maketitle
\begin{abstract}
This is the first part of a series of articles where we are going
to develop theory of valuations on manifolds generalizing the
classical theory of continuous valuations on convex subsets of an
affine space. In this article we still work only with linear
spaces. We introduce a space of smooth (non-translation invariant)
valuations on a linear space $V$. We present three descriptions of
this space. We describe the canonical multiplicative structure on
this space generalizing the results from \cite{alesker-poly}
obtained for polynomial valuations.
\end{abstract}
\setcounter{section}{-1}
\section{Introduction.}\setcounter{subsection}{1}
This is the first part of a series of articles where we are going
to develop theory of valuations on manifolds generalizing the
classical theory of continuous valuations on convex subsets of an
affine space. In this article we still work only with linear
spaces. In the subsequent parts of this series we are going to
generalize constructions of this article to arbitrary smooth
manifolds \cite{part2}, \cite{part3}. The case of a linear space
considered here will be useful for the general case for technical
reasons.

Let us remind some basic definitions. Let $V$ be a finite
dimensional real vector space, $n=\dim V$. Let $\ck(V)$ denote the
class of all convex compact subsets of $V$. Equipped with the
Hausdorff metric, the space $\ck(V)$ is a locally compact space.

\begin{definition}
a) A function $\phi :{\cal K}(V) \to \CC$ is called a valuation if
for any $K_1, \, K_2 \in {\cal K}(V)$ such that their union is
also convex one has
$$\phi(K_1 \cup K_2)= \phi(K_1) +\phi(K_2) -\phi(K_1 \cap K_2).$$

b) A valuation $\phi$ is called continuous if it is continuous
with respect to the Hausdorff metric on ${\cal K}(V)$.
\end{definition}

For the classical theory of valuations we refer to the surveys
McMullen-Schneider \cite{mcmullen-schneider} and McMullen
\cite{mcmullen-survey}. For the general background from convexity
we refer to Schneider \cite{schneider-book}.

In this article we introduce the space $SV(V)$ of {\itshape
smooth} valuations on $V$ (see Definition \ref{aff-1-2}). $SV(V)$
is a Fr\'echet space. We present three different descriptions of
this space. We describe a canonical structure of commutative
associative topological algebra with unit (where the unit is the
Euler characteristic).

Moreover the algebra $SV(V)$ has a canonical filtration by closed
subspaces
$$SV(V)=W_0\supset W_1\supset\dots\supset W_n$$
compatible with the product, namely $W_i\cdot W_j\subset W_{i+j}$.
Note that the subspace $W_n$ coincides with the space of smooth
densities on $V$. Moreover in Theorem \ref{mult-1-3} we prove that
there exists canonical isomorphism of the associated graded
algebra $gr_W SV(V):=\bigoplus_{i=0}^nW_i/W_{i+1}$ and the algebra
$C^\infty(V,Val^{sm}(V))$ of infinitely smooth functions on $V$
with values in the algebra $Val^{sm}(V)$ of smooth translation
invariant valuations.

Note also that the space $SV(V)$ contains polynomial smooth
valuations (studied by the author in \cite{alesker-poly}) as a
dense subspace. Thus the above results generalize results on
polynomial valuations from \cite{alesker-poly}.

The paper is organized as follows. In Section \ref{background} we
remind some necessary facts from the representation theory
(Subsection \ref{rep}) and the valuation theory (Subsection
\ref{val}).

In Section \ref{smooth-val} we introduce the main object of this
article, namely the space of smooth valuations (Definition
\ref{aff-1-2}).

In Section \ref{aff-3} we introduce the filtration $W_\bullet$ on
$SV(V)$. We study its basic properties. In particular we show in
Proposition \ref{aff-3-5} the isomorphism of Fr\'echet spaces
$$gr_W SV(V)\simeq C^\infty(V,Val^{sm}(V)).$$
This isomorphism is the first description of the space of smooth
valuations. In Corollary \ref{aff-3-9} we obtain the second one.

In Section \ref{mult} we introduce and study the canonical
multiplicative structure on $SV(V)$.

In Section \ref{intcc} we remind the construction of continuous
valuations using integration with respect to the normal cycle.
Then we show in Theorem \ref{cc-2-2} that all smooth valuations
are obtained in this way. This is the third promised description
of smooth valuations.

{\bf Acknowledgements.} I am grateful to J. Bernstein to numerous
very useful discussions. I express my gratitude to J. Fu for very
fruitful conversations, and in particular for his explanations of
the construction of valuations using the integration over the
normal cycle. I thank V.D. Milman for his interest to this work
and useful discussions, and B. Mityagin and V. Palamodov for
useful discussions.

\section{Background.}\label{background} In this section we remind some necessary
facts from representation theory (Subsection \ref{rep}) and theory
of valuations (Subsection \ref{val}). No result of this section is
new.

\subsection{Some representation theory.}\label{rep} \setcounter{theorem}{0}
\begin{definition}\label{rep-1} Let $\rho$ be a continuous representation of a Lie group $G$ in a
Fr\'echet space $F$. A vector $\xi \in F$ is called $G$-smooth if
the map $g\mapsto \rho(g)\xi$ is an infinitely differentiable map
from $G$ to $F$.
\end{definition}
It is well known (see e.g. \cite{wallach}, Section 1.6) that the
subset $F^{sm}$ of smooth vectors is a $G$-invariant linear
subspace dense in $F$. Moreover it has a natural topology of a
Fr\'echet space (which is stronger than the topology induced from
$F$), and the representation of $G$ in $F^{sm}$ is continuous.
Moreover all vectors in $F^{sm}$ are $G$-smooth.

Let $G$ be a real reductive group. Assume that $G$ can be imbedded
into the group $GL_N(\RR)$ for some $N$ as a closed subgroup
invariant under the transposition. Let us fix such an imbedding
$p:G\hookrightarrow GL_N(\RR)$. (In our applications $G$ will be
either $GL_n(\RR)$ or a direct product of several copies of
$GL_n(\RR)$.) Let us introduce a norm $|\cdot |$ on $G$ as
follows:

$$|g|:=\max\{||p(g)||,||p(g^{-1})||\}$$
where $||\cdot||$ denotes the usual operator norm in $\RR^N$.
\begin{definition}
Let $(\pi,G,F)$ be a smooth representation of $G$ in a Fr\'echet
space $F$ (namely $F^{sm}=F$). One says that this representation
has {\itshape moderate growth} if for each continuous semi-norm
$\lambda$ on $F$ there exists a continuous semi-norm $\nu_\lambda$
on $F$ and $d_{\lambda}\in \RR$ such that
$$\lambda(\pi(g)v)\leq ||g||^{d_\lambda}\nu_{\lambda}(v)$$
for all $g\in G,\, v\in F$.
\end{definition}

The proof of the next lemma can be found in \cite{wallach}, Lemmas
11.5.1 and 11.5.2.
\begin{lemma}\label{wallach}
(i) If $(\pi,G,H)$ is a continuous representation of $G$ in a
Banach space $H$ then $(\pi,G,H^{sm})$ has moderate growth.

(ii) Let $(\pi, G,V)$ be a representation of moderate growth. Let
$W$ be a closed $G$-invariant subspace of $V$. Then $W$ and $V/W$
have moderate growth.
\end{lemma}
The next lemma is obvious.
\begin{lemma}\label{obvious}
Let $G_1$ be a closed reductive subgroup of a reductive group $G$.
Assume that the image of $G_1$ in $GL_N(\RR)$ under the map
$p:G\hookrightarrow GL_N(\RR)$ is closed under the transposition.
Let $(\pi,G,F)$ has moderate growth. Then the restriction of this
representation to $G_1$ also has moderate growth.
\end{lemma}

Remind that a continuous Fr\'echet representation $(\rho,G,\cf)$
is said to have {\itshape finite length} if there exists a finite
filtration
$$0=F_0\subset F_1\subset \dots\subset F_m=F$$
by $G$-invariant closed subspaces such that $F_i/F_{i-1}$ is
irreducible, i.e. does not have proper closed $G$-invariant
subspaces.

A Fr\'echet representation $(\rho,G,F)$ of a real reductive group
$G$ is called {\itshape admissible} if its restriction to a
maximal compact subgroup $K$ of $G$ contains an isomorphism class
of any irreducible representation of $K$ with at most finite
multiplicity. (Remind that a maximal compact subgroup of
$GL_n(\RR)$ is the orthogonal group $O(n)$.)

\begin{theorem}[Casselman-Wallach, \cite{casselman}]\label{casselman-wallach}
Let $G$ be a real reductive group. Let $(\rho,G,F_1)$ and
$(\pi,G,F_2)$ be smooth representations of moderate growth in
Fr\'echet spaces $F_1, F_2$. Assume in addition that $F_2$ is
admissible of finite length. Then any continuous morphism of
$G$-modules $f:F_1\to F_2$ has closed image.
\end{theorem}
Let us also remind the classical L. Schwartz kernel theorem.
\begin{theorem}[L. Schwartz kernel theorem,
\cite{gelfand-vilenkin}]\label{schwartz}
 Let $X_1$and $X_2$ be compact smooth
manifolds. Let $\ce_1$ and $\ce_2$ be smooth finite dimensional
vector bundles over $X_1$ and $X_2$ respectively. Let $\cg$ be a
Fr\'echet space. Let
$$B:C^\infty(X_1,\ce_1)\times C^\infty(X_2,\ce_2)\to \cg$$
be a continuous bilinear map. Then there exists unique continuous
linear operator
$$b:C^\infty(X_1\times X_2,\ce_1\boxtimes \ce_2)\to \cg$$
such that $b(f_1\otimes f_2)=B(f_1,f_2)$ for any $f_i\in
C^\infty(X_i,\ce_i),\, i=1,2$.
\end{theorem}

The proof of L. Schwartz kernel theorem is based on the next
elementary and well known lemma which will be used in this
article.
\begin{lemma}\label{techn-lemma}
Let $X_1,X_2$ be two smooth manifolds such that $X_2$ is compact.
Let $\ce_1$ and $\ce_2$ be smooth finite dimensional vector
bundles over $X_1$ and $X_2$ respectively. Let $M\in \NN$ be an
integer. Let $G\subset X_1$ be a compact subset.

Then there exists a compact subset $\tilde G\subset X_1$
containing $G$, an integer $N\in \NN$, and a constant $C$ such
that for any $f\in C^\infty(X_1\times X_2,\ce_1\boxtimes \ce_2)$
there exists a presentation
$$f=\sum_{i=1}^\infty g_i\otimes h_i$$
such that $g_i\in C^\infty(X_1,\ce_1),\, h_i\in
C^\infty(X_2,\ce_2)$ and
$$\sum_{i=1}^\infty ||g_i||_{C^M(G)}||h_i||_{C^M(X_2)}\leq
C||f||_{C^N(\tilde G\times X_2)}.$$
\end{lemma}

For a Fr\'echet space $F$ and a smooth manifold $X$ let us denote
by $C^\infty(X,F)$ the Fr\'echet space of infinitely smooth
$F$-valued functions on $X$ with the topology of uniform
convergence with all derivatives on compact subsets of $X$. The
next proposition is well known but we do not have a reference.
\begin{proposition}\label{epi}
Let $G$ be a real reductive Lie group. Let $F_1,F_2$ be continuous
Fr\'echet $G$-modules. Let $\xi:F_1\to F_2$ be a continuous
morphism of $G$-modules. Assume that the assumptions of the
Casselman-Wallach theorem are satisfied, namely $F_1$ and $F_2$
are smooth and have moderate growth, and $F_2$ is admissible of
finite length. Assume moreover that $\xi$ is surjective.

Let $X$ be a smooth manifold. Consider the map
$$\hat \xi:C^\infty(X,F_1)\to C^\infty(X,F_2)$$
defined by $(\hat\xi (f))(x)=\xi(f(x))$ for any $x\in X$.

Then $\hat\xi$ is surjective.
\end{proposition}
{\bf Proof.} First  let us prove this proposition under assumption
that $X=\PP^n$. Consider the natural action of the group
$GL_{n+1}(\RR)$ on $\PP^n$. Then the representation of the
reductive group $G\times GL_{n+1}(\RR)$ in the spaces
$C^\infty(X,F_i),\, \, i=1,2$ is smooth and of moderate growth.
Moreover $C^\infty(X,F_2)$ is an admissible $(G\times
GL_{n+1}(\RR))$-module of finite length. Hence by the
Casselman-Wallach theorem it is enough to show that $\hat \xi$ has
 dense image. But since $\xi$ is surjective, the image of $\hat\xi$ contains finite
 linear combinations of elements of the form $f\otimes v$ with
 $f\in C^\infty(\PP^n),\, v\in F_2$. Clearly such linear
 combinations are dense in $C^\infty(X,F_2)$.

Let us return to the case of a general manifold $X$. Let us denote
by $\co_X$ the sheaf of infinitely smooth functions on $X$. Let us
consider the sheaves $\cf_i,\, i=1,2,$ on $X$ defined by
$$\cf_i(U)=C^\infty(U,F_i)$$
for any open subset $U\subset X$. Then $\xi$ induces a morphism of
$\co_X$-modules (which will be also denoted by $\hat\xi$)
$$\hat\xi:\cf_1\to \cf_2$$
which is obviously defined. Let us show that $\hat\xi$ is an
epimorphism of sheaves. Since this statement is local, and all
smooth manifolds of given dimension are locally diffeomorphic, we
may assume again that $X=\PP^n$. Fix a point $x\in \PP^n$. Let
$B_1\subset B_2\subset \PP^n$ be two open balls such that $x\in
B_1$ and the closure of $B_1$ is contained in $B_2$. Let us fix a
function $\gamma\in C^\infty(\PP^n)$ such that
$\gamma|_{B_1}\equiv 1,\, \gamma|_{\PP^n\backslash B_2}\equiv 0$.
Let $\phi\in H^0(B_2,\cf_2)$. Set $\psi:=\gamma \cdot \phi$. Then
$\psi$ extend by zero to a section from $H^0(\PP^n,\cf_2)=F_2$.
This section will be denoted again by $\psi$. By the previous case
there exists $\chi\in F_1= H^0(\PP^n,\cf_1)$ such that
$\hat\xi(\chi)=\psi$. Restricting the last identity to $B_1$ we
conclude that $\hat\xi$ is an epimorphism of sheaves.

For a general $X$, let us denote by $\ck:=Ker \hat\xi$. Thus $\ck$
is an $\co_X$-module. It is well known that every $\co_X$-module
$\ck$ is acyclic, i.e. $H^i(X,\ck)=0$ for $i>0$ (see e.g.
\cite{godement}).

We have a short exact sequence of sheaves
$$0\to \ck\to \cf_1\to \cf_2\to 0.$$
From the long exact sequence we get
$$H^0(X,\cf_1)\to H^0(X,\cf_2)\to H^1(X,\ck)=0.$$
The result is proved. \qed

\begin{remark}
Recently we were informed by B. Mityagin and V. Palamodov that the
following general fact, which is also sufficient for our purposes
instead of Proposition \ref{epi}, is true. Let $\xi:F_1\to F_2$ be
a surjective continuous linear map of nuclear Fr\'echet spaces.
Let $X$ be a smooth manifold. Then $\hat\xi:C^\infty(X,F_1)\to
C^\infty(X,F_2)$ is surjective.
\end{remark}

\subsection{Some valuation theory.}\label{val}\setcounter{theorem}{0} Let us
remind few basic facts from the theory of translation invariant
continuous valuations. For a real vector space $V$ of finite
dimension $n$ let us denote by $Val(V)$ the space of translation
invariant valuations on $\ck(V)$ continuous with respect to the
Hausdorff metric on $\ck(V)$. Equipped with the topology of
uniform convergence on compact subsets of $\ck(V)$ the space
$Val(V)$ becomes a Banach space (see e.g. Lemma A.4 in
\cite{alesker-poly}).

\begin{definition}
A valuation $\phi$ is called homogeneous of degree $k$ (or just
$k$-homogeneous) if for every convex compact set $K$ and for every
scalar $\lam
>0$
$$\phi(\lam K)=\lam ^k \phi (K).$$
\end{definition}
Let us denote by $Val _k(V)$ the space of translation invariant
continuous valuations homogeneous of degree $k$.
\begin{theorem}[McMullen \cite{mcmullen-euler}]
$$Val(V)=\bigoplus_{k=0}^{n} Val_k(V),$$
where $n=\dim V$.
\end{theorem}
Note in particular that the degree of homogeneity is an integer
between 0 and $n=\dim V$. It is known that $Val_0(V)$ is
one-dimensional and is spanned by the Euler characteristic $\chi$,
and $Val_n(V)$ is also one-dimensional and is spanned by a
Lebesgue measure \cite{hadwiger-book}. The space $Val_n(V)$ is
also denoted by $|\wedge V^*|$ (the space of complex valued
Lebesgue measures on $V$). One has further decomposition with
respect to parity:
$$Val_k(V)=Val_k^{ev}(V)\oplus Val_k^{odd}(V),$$
where $Val_k^{ev}(V)$ is the subspace of even valuations ($\phi$
is called even if $\phi(-K)=\phi(K)$ for every $K\in {\cal
K}(V)$), and $Val_k^{odd}(V)$ is the subspace of odd valuations
($\phi$ is called odd if $\phi (-K)=-\phi(K)$ for every $K\in
{\cal K}(V)$). The Irreducibility Theorem is as follows.
\begin{theorem}[Irreducibility Theorem
\cite{alesker-gafa}]\label{irr} The natural representation of the
group $GL(V)$ on each space $Val_k^{ev}(V)$ and $Val_k^{odd}(V)$
is irreducible for any $k=0,1,\dots,n$.
\end{theorem}
In this theorem, by the natural representation one means the
action of $g\in GL(V)$ on $\phi\in Val(V)$ as
$(g\phi)(K)=\phi(g^{-1}K)$ for every $K\in \ck(V)$. The subspace
of smooth valuations with respect to this action in sense of
Definition \ref{rep-1} is denoted by $Val^{sm}(V)$.
\begin{remark}
The representation $Val(V)$ of $GL(V)$ is an admissible
representation. Indeed, it was show in \cite{alesker-advmath} that
$Val_k^{ev/odd}$ can be $GL(V)$-equivariantly imbedded into the
space of continuous sections of a $GL(V)$-equivariant finite
dimensional vector bundle $\ce^{ev/odd}$ over the projective space
$\PP_+(V^*)$. Let us show that this representation must be
admissible.  Let us fix on $V^*$ a Euclidean metric. Let us fix a
point $l_0\in \PP(V^*)$. Let $H$ denote the stabilizer of $l_0$ in
$O(n)$. Then $\PP_+(V^*) \simeq O(n)/H$. Let $q\colon O(n)\to
\PP(V^*)$ be the surjection $g\mapsto g(l_0)$. Let
$\cl:=q^*\ce^{ev/odd}$. Then $\cl$ is $O(n)$-equivariant vector
bundle over $O(n)$, hence $\cl$ is $O(n)$-equivariantly trivial.
Note that we have the $O(n)$-equivariant imbedding
$C(\PP(V^*))\hookrightarrow C(O(n),\cl)$. Hence it is enough to
check that $C(O(n),\cl)$ contains each irreducible representation
of $O(n)$ with at most finite multiplicity. This follows from the
well known fact that for any compact group $K$ the space of
functions $C(K)$ contains each irreducible representation $\pi$ of
$K$ with multiplicity $\dim \pi$ (which is necessarily finite).
\end{remark}

\section{The space of smooth valuations}\label{smooth-val}\setcounter{subsection}{0}
\subsection{Some definitions.}\label{aff-1} Let $V$ be an
$n$-dimensional real vector space. Let us denote by $CV(V)$ the
space of continuous valuation on $V$. Equipped with the topology
of uniform convergence on compact subsets of $\ck(V)$, $CV(V)$
becomes a Fr\'echet space. Let $\qvv$ denote the space of
continuous valuations on $V$ which satisfy the following
additional property:
\newline
the map given by $K\mapsto \phi( tK +x)$ is a continuous map $\ck
(V) \to C^{n}([0,1]\times V)$. Let us call such valuations
{\itshape quasi-smooth}.

In the space $\qvv$ we have the natural linear topology defined as
follows. Fix a compact subset $G\subset V$. Define a seminorm on
$\qvv$
$$||\phi||_G:=sup \{||\phi(tK+x)||_{C^{n}([0,1]\times G)}|\, K\subset G\}.$$
Note that the seminorm $||\cdot||_G$ is finite. One easily checks
the following claim.
\begin{claim}\label{aff-1-1}
Equipped with the topology defined by this sequence of seminorms
the space $\qvv$ is a Fr\'echet space.
\end{claim}

Note also that the natural representation of the group $Aff(V)$ of
affine transformations of $V$ in the space $\qvv$ is continuous.
We will denote by $\svv$ the subspace of $Aff(V)$-smooth vectors
in $\qvv$. It is a Fr\'echet space.

\begin{definition}\label{aff-1-2}
Elements of $SV(V)$ will be called {\itshape smooth} valuations on
$V$.
\end{definition}

\subsection{Main examples.}\label{aff-2}
\setcounter{theorem}{0}
 Let $V$ be a real vector space of
dimension $n$. Let us denote by $\pl$ the manifold of oriented
lines passing through the origin in $V^*$. Let $L$ denote the line
bundle over $\pl$ such that its fiber over an oriented line $l$
consists of linear functionals on $l$. Let $\mov$ denote the line
bundle of densities over $V$. Let $p:V\times \pl \to V$ be the
projection. For any integer $k$, $0\leq k\leq n$, we are going to
construct a natural map
$$\Theta_k: \bigoplus_{j=0}^kC^{\infty}(V\times (\pl)^j, \mov\boxtimes L^{\boxtimes
j})\to \sv(V).$$

First let us remind some results from \cite{alesker-int}.  Let
$\bar K=(K_1 , K_2 ,\dots, K_s )$ be an $s$-tuple of compact
convex subsets of $V $. Let $r\in \NN\cup \{\infty\}$. For any
$\mu \in C^{r}(V,\mov)$ consider the function $M_{\bar K} \mu \,
:\RR^{s}_{+} \to \CC \, ,\mbox{where }\RR_{+}^{s}=
\{(\lam_1,\dots,\lam_s) \,| \,\lam_j\geq0 \mbox{ for all } j\} $
defined by
$$ ( M_{\bar K} \mu )
  (\lam _1,\dots,\lam _s)= \mu(\sum_{j=1}^{s} \lam _j K_j) .  $$
\begin{theorem}[\cite{alesker-int}]\label{aff-2-1}
(1) $M_{\bar K}\mu \in C^{r}(\RR^s_+)$ and $M_{\bar K}$ is a
continuous operator from $C^{r}(V,\mov)$ to $C^{r}(\RR^s_+)$.

(2) Assume that  a sequence $\mu^{(m)}$ converges to $\mu$ in
$C^{r}(V,\mov)$. Let $ K_{j}^{(m)}, \, K_{j},\,
j=1,\dots,s,\,m\in\NN $, be convex compact sets in $ V $, and for
every $ j=1,\dots ,s $ $K_{j}^{(m)} \to K_{j} $ in the Hausdorff
metric as $m\to \infty$. Then $ M_{\bar K^{(m)} } \mu^{(m)} \to
M_{\bar K } \mu $ in $C^{r}(\RR^s_+)$.
\end{theorem}

Before we define the map $\Theta _k$  let us make more remarks.
Fix $s$, $0\leq s\leq k$.  Let us fix $\mu \in C^{\infty}(V,\mov)$
and $A_1,\dots, A_s\in \ck(V)$ being strictly convex with smooth
boundaries . Let us define
$$(\Theta_{s}'(\mu;A_1,\dots,A_s))(K):=\frac{\pt^s}{\pt
\lam_1\dots\pt\lam_s}\big |_{0}\mu(K+\sum_{j=1}^s \lam_jA_j).$$

 Theorem \ref{aff-2-1} implies that $\Theta_{s}'(\mu;A_1,\dots,A_s)\in
\sv(V)$. It is clear that $\Theta_{s}'$ is Minkowski additive with
respect to each $A_j$. Namely, say for $j=1$, one has
$\Theta_{s}'(\mu;a A_1'+b
A_1'',A_2,\dots,A_s)=a\Theta_{s}'(\mu;A_1',A_2,\dots,A_s)+
b\Theta_{s}'(\mu;A_1'',A_2,\dots,A_s)$ for $a,\, b \geq 0$.

Remind that for any $A\in \ck(V)$ one defines the supporting
functional $h_A(y):=\sup_{x\in A}y(x)$ for any $y\in V^*$. Thus
$h_A\in C(\PP_+(V^*),L)$. Moreover it is well known (and easy to
see) that $A_N\to A$ in the Hausdorff metric if and only if
$h_{A_N}\to h_A$ in $C(\PP_+(V^*),L)$. Also any section $F\in
C^2(\PP_+(V^*),L)$ can be presented as a difference $F=G-H$ where
$F,\,H\in C^2(\PP_+(V^*),L)$ are supporting functionals of some
convex compact sets and $\max\{||G||_2,||H||_2\}\leq c||F||_2$
where a constant $c$ is independent of $F$. (Indeed one can choose
$G=F+R\cdot h_D,\, H=R\cdot h_D$ where $D$ is the unit Euclidean
ball, and $R$ is a large enough constant depending on $||F||_2$.)
Hence we can uniquely extend $\Theta_{s}'$ to a multilinear
continuous map (which we will denote by the same letter):
$$\Theta_{s}':C^{\infty}(V,\mov) \times (C^{\infty}(\PP_+(V^*),L))^s\to \sv(V).$$
Theorem \ref{aff-2-1} implies that $\Theta_s'$ depends
continuously on each argument. By the L. Schwartz kernel theorem
it follows that this map gives rise to a continuous linear map
$$\Theta_{s}': C^{\infty}(V\times \PP_+(V^*)^s,\mov\boxtimes L^{\boxtimes
s})\to \sv(V).$$ Now let us define the map
$$\Theta_k:=\bigoplus_{i=0}^{k} \Theta _i'.$$


\section{Filtrations.}\label{aff-3}\setcounter{subsection}{1}
\setcounter{theorem}{0}
 Let us define a decreasing filtration on
$\svv$. Set
 $$W_i:= \{\phi \in \svv|\, \frac{d^k}{dt^k}
  \phi(tK+x)\big |_{t=0}=0 \,
 \forall k<i,\, K\in \ck(V),\, x\in V\}.$$
 It is clear that $W_i$ are $Aff(V)$-invariant closed subspaces of
 $\svv$.
Obviously $\svv =W_0\supset W_1\supset \dots .$
\begin{proposition}\label{aff-3-1}
$$W_{n+1}=0.$$
\end{proposition}
\def\tj{T_{j_1 \dots j_{l-1}}}
{\bf Proof.} Let $\phi\in W_{n+1}$. We want to show that $\phi$
vanishes. Let us prove it by induction in $n=\dim V$. For $n=0$
the statement is clear. Let us assume that the statement holds for
$n-1$. Then this implies that $\phi$ is a simple valuation, i.e.
it vanishes on convex sets of dimension less than $n$. It is
sufficient to show that $\phi$ vanishes on polytopes. Since every
polytope can be dissected into simplices it is sufficient to prove
that $\phi$ vanishes on simplices. Let $\Delta$ be  a simplex.
Choosing an appropriate coordinate system we may assume that it
has the form
$$\Delta=\{(x_1,\dots, x_n)|\, 0\leq x_1\leq \dots \leq x_n\leq
1\}.$$ Set for $ 1 \leq i \leq j \leq n $, $ T_{i,j} := $
$$ \left\{ ( x_{1}, \dots , x_{n} ) \, | \, 0 \leq x_{i} \leq  x_{i+1} \leq \dots \leq x_{j} \leq 1
\mbox{ and }  x_{l} =0
 \mbox { for } l< i \mbox
{ and } l> j \right\} . $$ For a sequence $ 0< j_{1} < \dots <
j_{l-1 } <n $, let us denote (as in \cite{mcmullen-euler} )
 $$ \tj :=
T_{0 j_{1} } + \dots + T_{j_{l-1}, \, n }. $$ First note that any
point $z\in \Delta$ has the form $ z=(z_{i} ) _{i=1} ^{n} $, where
\begin{equation} z_{1} = \dots =z_{j_{1} } < z_{ j_{1} +1 } =
\dots= z_{ j_{2} } <
 \dots
< z_{j_{l-1} +1 } = \dots = z_{j_{l} } \leq 1  ,
\end{equation}
and $ j_{l} = n $.

\def\rj{R_{j_1\dots j_{l-1}}}
For a sequence $ 0< j_{1} < \dots < j_{l-1 } <n $ let us also
define $$R_{j_1\dots j_{l-1}}(N):=\{z\in \frac{1}{N}\ZZ ^n\cap
\Delta |\, z \mbox{ satisfies } (1)\}.$$ Then since $\phi$ is a
simple valuation one has
\begin{equation} \phi(\Delta)=\sum
_{0<j_1<\dots < l_{l-1}<n}\left(\sum _{z\in \rj
(N)}\phi(z+\frac{1}{N}\tj)\right).
\end{equation}

Since $\phi\in W_{n+1}$, for any $\eps>0$ there exists $N(\eps)$
such that for all $N>N(\eps)$ and for all $z\in \Delta$ one has
$|\phi(z+\frac{1}{N}\tj)|<\eps N^{-n}$. However it is easy to see
that $|\frac{1}{N}\ZZ^n\cap \Delta|\leq CN^n$ where $C$ is a
constant independent of $N$. From the equation $(2)$ we get an
estimate $|\phi(\Delta)|\leq C' \eps$ where $C'$ is a constant
depending on $n$ only. Hence $\phi (\Delta)=0$. Hence $\phi \equiv
0$. \qed


\begin{proposition}\label{aff-3-2}
$W_n$ coincides with the space of smooth densities on $V$.
\end{proposition}
{\bf Proof.}
 Obviously smooth densities are contained in
$W_n$. Now let us fix $\phi\in W_n$. Let $\tilde
\phi(K,x)=\frac{d^n}{dt^n}\big |_{t=0} \phi(tK+x)$. Then for a
fixed $x\in V$, $\tilde\phi(\cdot, x)$ is a translation invariant
continuous valuation homogeneous of degree $n$. To check it, fix
an arbitrary compact subset $G\subset V$.  We have
$$\phi(tK+x)=\frac{t^n}{n!} \tilde\phi(K,x)+o(t^k)$$
uniformly on $K\subset G,\, t\in G$. Also $\tilde\phi(K,x)$
depends smoothly on $x$ when $K$ is fixed. Then
$\tilde\phi(t(K+a)+x)=\frac{t^n}{n!}\tilde\phi(K,ta+x)+o(t^n)=
\frac{t^n}{n!}\tilde\phi(K,x)+o(t^n).$
Hence $\tilde\phi(K,x)$ is translation invariant in $K$ when $x$
is fixed.

By a result due to Hadwiger \cite{hadwiger-book}
$\tilde\phi(\cdot,x)$ must be a Lebesgue measure. Also it depends
smoothly on $x\in V$. Subtracting from $\phi$ an appropriate
density and using the fact that $W_{n+1}=0$ by Proposition
\ref{aff-3-1}, we deduce the result. \qed

\begin{example}\label{aff-3-3}
(1) Let us remind the definition of a {\itshape polynomial
valuation} introduced by Khovanskii and Pukhlikov
\cite{khovanskii-pukhlikov1}. A valuation $\phi$ is called
polynomial of degree at most $d$ if for any $K\in \ck(V)$ the
function $V\to \CC$ given by $x\mapsto \phi(K+x)$ is a polynomial
of degree at most $d$. In \cite{khovanskii-pukhlikov1} it was
shown that if $\phi$ is a continuous polynomial valuation of
degree $d$ then for any $K_1,\dots, K_s \in \ck(V)$ the function
$\phi(\sum_j \lam_j K_j)$ is a polynomial in $\lam_j\geq 0$ of
degree at most $d+n$. It follows that $\phi \in \qvv$.

(2) Let $\mu$ be a smooth density. It was shown in
\cite{alesker-int} that the map
$$(K_1,\dots, K_s; \lam_1,\dots, \lam_s)\mapsto \mu(\sum_{j=1}^s
\lam _jK_j)$$ defines a continuous map $\ck(V)^s\to
C^{\infty}(\RR_{+}^s)$. It follows that for any fixed $A_1,\dots ,
A_s\in \ck(V)^s$ the map
$$K\mapsto \frac{\pt^s}{\pt \lam _1 \dots \pt \lam _s}\big |_{0}
\mu(K+\sum _j \lam_j A_j)$$ defines a valuation from $QV(V)$.
\end{example}

Remind that in Section \ref{aff-2} we have defined a map
$$\Theta_k:\bigoplus_{i=0}^k C^{\infty}(V\times \pl ^i,
\mov \boxtimes L^{\boxtimes i})\to \sv (V).$$

\begin{proposition}\label{aff-3-4}
The image of $\Theta _k$ is contained in $W_{n-k}$.
\end{proposition}
{\bf Proof.} Indeed $\mu(rK+x+\sum_j\lam_j A_j)=
O((\sqrt{r^2+\sum_j\lam_j^2})^n).$ The result follows from the
construction of $\Theta _k$. \qed

In Corollary \ref{aff-3-9} we will prove that in fact the image of
$\Theta_k$ coincides with $W_{n-k}$. Let us denote by $Val(TV)$
the (infinite dimensional) vector bundle over $V$ whose fiber over
$x\in V$ is equal to the space of translation invariant
$GL(T_xV)$-smooth valuations on the tangent space $T_xV$.
Similarly we can define the vector bundle $Val_k(TV)$ of
$k$-homogeneous smooth translation invariant valuations. Clearly
$C^\infty(V,Val_k(TV))=C^\infty(V,Val_k^{sm}(V))$ where the last
space denotes the space of $C^\infty$-smooth functions on $V$ with
values in the Fr\'echet space $Val_k^{sm}(V)$ of $k$-homogeneous
translation invariant smooth valuations.

Let us define a map
$$\Lambda_k: W_k\to C^{\infty}(V,Val_k(TV))$$ by
$\Lambda_k(\phi):=[ K\mapsto [x\mapsto
\frac{1}{k!}\frac{d^k}{dt^k}\big |_{t=0} \phi(tK+x)]]$.
\begin{proposition}\label{aff-3-5}
 (i) $\Lambda_k: W_k\to C^{\infty}(V,Val_k(TV))$ is an epimorphism.

 (ii) $Ker \Lambda_k=W_{k+1}$.

 (iii) $W_k/W_{k+1}$ is isomorphic to
$C^{\infty}(V,Val_k(TV))$.
\end{proposition}
{\bf Proof.} Clearly (iii) follows from (i) and (ii). Part (ii) is
obvious from the definitions.

 Let us check next that $K\mapsto
\frac{1}{k!}\frac{d^k}{dt^k}\big |_{t=0} \phi(tK+x)$ is a
translation invariant continuous valuation for any $x\in V$, and
it depends smoothly on $x$. The only thing one should check is the
translation invariance. Let us denote
$\psi(K,x):=\frac{1}{k!}\frac{d^k}{dt^k}\big |_{t=0} \phi(tK+x)$.
Then for any fixed compact subset $G\subset V$ we have
$$\phi(tK+x)=t^k \psi(K,x)+o(t^k)$$
uniformly in $K\subset G,\, x\in G$. Also $\psi(K,x)$ depends
smoothly on $x$. Then
$\psi(t(K+a)+x)=t^k\psi(K,ta+x)+o(t^k)=t^k\psi(K,x)+o(t^k)$. This
proves the translation invariance of the limit.

It remains to prove surjectivity of $\Lambda_k$. We will need the
following lemma.
\begin{lemma}\label{aff-3-6}
The map $$\Xi_k:=\Lambda_k\circ \Theta_{n-k}': C^{\infty}(V\times
\pl^{n-k},\mov \boxtimes L^{\boxtimes (n-k)})\to C^{\infty}(V,
Val_k(TV))$$ is an epimorphism.
\end{lemma}
Obviously Proposition \ref{aff-3-5}(i) follows from Lemma
\ref{aff-3-6}.

{\bf Proof} of Lemma \ref{aff-3-6}. Let us describe $\Xi_k$
explicitly. Let $\gamma=\mu\otimes h_{A_1}\otimes \dots \otimes
h_{A_{n-k}}$ where $h_{A_i}$ is the supporting functional of a set
$A_i\in \ck (V)$, $\mu\in C^{\infty}(V,|\ome_V|)$. Let $\mu=F(y)
dy$ where $dy$ is a Lebesgue measure. Then
$$(\Theta_{n-k}'\gamma)(K)=\frac{\pt^{n-k}}{\pt\lam_1\dots \pt
\lam_{n-k}}\big |_0 \int_{K+\sum _i\lam _i A_i} F(y) dy.$$ Hence
\begin{eqnarray*}
(\Xi_k \gamma)(K)=\lim_{r\to +0} \frac{1}{r^k}
\frac{\pt^{n-k}}{\pt \lam_1\dots \pt \lam_{n-k}}
\big|_0\int_{rK+x+\sum _i\lam _i A_i} F(y) dy=\\
\lim_{r\to +0}\frac{1}{r^k} \frac{\pt^{n-k}}{\pt \lam_1\dots \pt
\lam_{n-k}} \big|_0
\left( F(x) vol( rK+\sum _i\lam _i A_i)+
o\left(\left(\sqrt{r^2+\sum_i\lam_i^2}\right)^n\right)\right)=\\
\lim_{r\to +0}\frac{1}{r^k} \frac{\pt^{n-k}}{\pt \lam_1\dots \pt
\lam_{n-k}} \big|_0 F(x) vol( rK+\sum _i\lam _i A_i)=\\
F(x)\frac{\pt^{n-k}}{\pt \lam_1\dots \pt \lam_{n-k}} \big|_0
vol(K+\sum _i\lam _i A_i).
\end{eqnarray*}

This computation shows that $\Xi_k$ is a morphism of
$C^{\infty}(V)$-modules.
Let us denote
$$F_1:=C^\infty(\PP_+(V^*)^{n-k},|\wedge^n(V^*)|\otimes
L^{\boxtimes (n-k)}).$$ It is easy to see that
$$C^\infty(V\times\PP_+(V^*)^{n-k},|\ome_V|\boxtimes
L^{\boxtimes(n-k)})=C^\infty(V,F_1).$$
\def\ttk{\tilde\Theta_{n-k}}
Moreover $$(\Xi_k f)(v)=\ttk(f(v)) \, \, \forall v\in V,f\in F_1$$
where $$\ttk:F_1\to Val_k^{sm}(V)$$ is the restriction of
$\Theta_{n-k}'$ to the space $F_1$ which coincides with the
subspace of $C^\infty(V\times\PP_+(V^*)^{n-k},|\ome_V|\boxtimes
L^{\boxtimes(n-k)})$ consisting of elements invariant with respect
to translations to vectors from $V$.

The map $\ttk$ commutes with the natural action of the group
$GL(V)$. Hence $\ttk$ is onto by Irreducibility Theorem \ref{irr}
and Casselman-Wallach Theorem \ref{casselman-wallach}. Proposition
\ref{epi} implies that $\Xi_k$ is onto as well. Thus Lemma
\ref{aff-3-6} is proved. \qed


\begin{corollary}\label{aff-3-9}
The image of the map $\Theta_k:\bigoplus_{i=0}^{k}
C^{\infty}(V\times \pl ^i, \mov \boxtimes L^{\boxtimes i})\to \sv
(V)$ is equal to $W_{n-k}$.
\end{corollary}
{\bf Proof.} By Proposition \ref{aff-3-4} the image of the map
$\Theta_k$ is contained in $W_{n-k}$. Let us prove the opposite
inclusion by the induction in $k$. For $k=0$ this is just
Proposition \ref{aff-3-2}. Let us assume that
$Im(\Theta_{k'})=W_{n-k'}$ for $k'<k$. Let $\phi \in W_{n-k}$. By
Proposition \ref{aff-3-5} there exists $\psi\in Im (\Theta_k)$
such that $\Lambda_{n-k}(\phi)=\Lambda_{n-k}(\psi)$. It follows
that $\phi-\psi \in W_{n-k+1}$. Applying the induction assumption
to this valuation we obtain the result. \qed

\begin{corollary}\label{aff-3-10}
Polynomial valuations from $W_k$ are dense in $W_k$. Polynomial
valuations from $QV(V)$ are dense in $QV(V)$.
\end{corollary}
{\bf Proof.} First notice that the second statement follows from
the first one. Indeed it is true since $W_0=\sv(V)$ is dense in
$QV(V)$. The first statement follows from Corollary \ref{aff-3-9}
and the obvious fact the image under $\Theta_k$ of any element of
$C^{\infty}(V\times \pl ^i, \mov \boxtimes L^{\boxtimes i})$ which
is polynomial with respect to translations in $V$, is  a
polynomial valuation. \qed

From this corollary we immediately get
\begin{corollary}\label{aff-3-11}
Let $G$ be a compact subgroup of $GL(V)$. Then $G$-invariant
polynomial valuations are dense in the space of $G$-invariant
quasi-smooth valuations.
\end{corollary}


Let us now introduce another decreasing filtration on $\sv(V)$.
Set
$$\gamma_i:=\{\phi \in \sv(V)|\, \phi(K)=0 \mbox{ if } \dim
K<i\}.$$ Clearly $\sv(V)=\gamma_0\supset \gamma_1\supset \dots
\supset \gamma_n \supset \gamma_{n+1}=0$.
\begin{theorem}\label{aff-3-12}
(i) $W_1=\gamma_1$.

(ii) $\gamma_{j+1}\subset W_j\subset \gamma_j$ for any $j$.
\end{theorem}
{\bf Proof.} (i) First note that for any $j$ we have $W_j\subset
\gamma_j$. This follows from Proposition \ref{aff-3-2} applied for
$(j-1)$-dimensional subsets. Let us prove that $\gamma_1\subset
W_1$. Let $\phi \in \gamma_1$, i.e. $\phi$ vanishes on points.
Hence for any $K\in \ck(V)$ the function $[0,1]\to C^n(V)$ given
by $t\mapsto [x\mapsto \phi(tK+x)]$ vanishes at $t=0$. Hence $\phi
\in W_1$ by the definition of $W_1$.

(ii) We have proven the second inclusion in (ii). Let us prove the
first one, namely $\gamma_{j+1}\subset W_j$. Assume that this is
not true. Then there exists $l<j$ and $\phi\in \gamma_{j+1}\cap
W_l$ such that $\phi \not \in W_{l+1}$. Set $\tilde
\phi:=\Lambda_l(\phi)\in C^\infty(V,Val_l^{sm}(TV))$. Then
$\tilde\phi \ne 0$ by Proposition \ref{aff-3-5}(ii). From the
construction of $\Lambda_l$ it follows that at each point $x\in V$
one has $\tilde \phi _x \in Val^{sm}_l(T_xV)\cap \gamma _{j+1}$.
But the last intersection vanishes; for translation invariant
valuations this was proved in \cite{alesker-poly} (see the
beginning of Section 3 in \cite{alesker-poly}). Thus we get a
contradiction. \qed

\section{The multiplicative structure.}\label{mult}\setcounter{subsection}{1}
\setcounter{theorem}{0}
 In this section we construct a
canonical multiplicative structure on $\sv(V)$. $SV(V)$ will
become a commutative associative algebra with unit (where the unit
is the Euler characteristic).

First of all we will construct for any linear spaces $X$ and $Y$
the exterior product
$$\sv(X)\times \sv(Y)\to QV(X\times Y)$$
which is a bilinear continuous map.

First let us introduce some notation. Let us denote for brevity
$$F_X:=\bigoplus_{k=0}^{\dim X}C^{\infty}(X\times
\PP_+(X^*)^k,|\ome_X|\boxtimes L^{\boxtimes k}),$$
$$F_Y:=\bigoplus _{l=0}^{\dim Y}C^{\infty}(Y\times
\PP_+(Y^*)^l,|\ome_Y|\boxtimes L^{\boxtimes l}).$$ First for any
$k\leq \dim X,\, l\leq \dim Y$ let us define a multilinear map
$$M: C^{\infty}(X,|\ome_X|)\times C^{\infty}(\PP_+(X^*),L)^k \times
C^{\infty}(Y,|\ome_Y|)\times C^{\infty}(\PP_+(Y^*),L)^l\to
QV(X\times Y).$$ Let $\mu\in C^{\infty}(X,|\ome_X|),\, \nu\in
C^{\infty}(Y,|\ome_Y|),\, \xi_i\in C^{\infty}(\PP_+(X^*),L),\,
\eta_j\in C^{\infty}(\PP_+(Y^*),L)$ where $i=1,\dots, k,\,
j=1,\dots l$. First let us assume that $\xi_i=h_{A_i}$ is a
supporting functional of a convex set $A_i\in \ck(X)$, and
$\eta_j=h_{B_j}$ is a supporting functional of $B_j\in \ck(Y)$.
Let us define
$$M(\mu,\xi_1,\dots, \xi_k;\nu,\eta_1,\dots,
\eta_l)(K)=$$ $$\frac{\pt^k}{\pt \lam_1\dots \pt
\lam_k}\frac{\pt^l}{\pt \theta_1\dots \pt \theta_l}\big |_0
(\mu\boxtimes \nu)(K+\sum_{i=1}^k\lam _i (A_i\times 0)+
\sum_{j=1}^l \theta_j(0\times B_j)).$$ It is clear that the right
hand side is Minkowski additive with respect to $A_i$ and $B_j$.
Hence using the same argument as in the construction of
$\Theta'_i$ (in Section \ref{aff-2}) and using Theorem
\ref{aff-2-1} we extend $M$ to a continuous multilinear functional
defined for all $\xi_i\in C^{\infty}(\PP_+(X^*),L),\, \eta_j\in
C^{\infty}(\PP_+(Y^*),L)$.

Hence by the L. Schwartz kernel theorem we get a bilinear
continuous map
$$M: F_X\times F_Y\to QV(X\times Y).$$
Remind that we have canonical surjections
$$\Theta_X:F_X\to \sv(X),\, \Theta_Y:F_Y\to \sv(Y)$$
where now we use the subscript to emphasize dependence on the
space.
\begin{lemma}\label{mult-1-1}
The bilinear map $M: F_X\times F_Y\to QV(X\times Y)$ admits a
unique factorization to a continuous bilinear map
$$M':\sv(X)\times\sv(Y)\to QV(X\times Y)$$
such that $M=M'\circ (\Theta_X\times \Theta_Y).$
\end{lemma}
{\bf Proof.} The uniqueness of such a factorization is obvious due
to the surjectivity of  $\Theta_X,\, \Theta_Y$. Let us prove
existence. Let us fix $f\in Ker \Theta_X$. It is enough to show
that $M(f,g)=0$ for any $g\in F_Y$. It is enough to assume that
$g=\nu\otimes\eta_1\otimes \dots \otimes \eta_p$ where $\nu\in
C^{\infty}(Y,|\ome_Y|),\, \eta_i\in C^{\infty}(\PP_+(Y^*),L)$, and
moreover $\eta_i=h_{B_i}$ where $B_i\in \ck(Y)$.

Let us prove that for any $w\in F_X$ and $K\in \ck(X\times Y)$ one
has \begin{equation}\label{diff} M(w,g)(K)=\frac{\pt ^p}{\pt
\theta_1\dots \pt\theta_p}\big |_0\int _{y\in Y}
\Theta_X(w)\left(\left(K+\sum_{j=1}^p\theta _j(0\times
B_j)\right)\cap (X\times\{y\})\right)d\nu(y).\end{equation} Note
that this identity implies Lemma \ref{mult-1-1}.

Let us check the identity (3).
First let us check it for $w=\mu\otimes \xi_1\otimes \dots \otimes
\xi_k$ where $\xi_i=h_{A_i}$, $A_i\in \ck(X)$. For such $w$ using
Theorem \ref{aff-2-1} one obtains
$$M(w,g)=\frac{\pt ^k}{\pt \lam_1\dots
\pt\lam_k}\frac{\pt ^p}{\pt \theta_1\dots \pt\theta_p}\big |_0
(\mu\boxtimes\nu)(K+\sum_{i=1}^k\lam_i(A_i\times 0)
+\sum_{j=1}^p\theta_j (0\times B_j))=$$
$$\frac{\pt ^p}{\pt \theta_1\dots \pt\theta_p}\frac{\pt ^k}{\pt \lam_1\dots
\pt\lam_k}\big |_0\int_{y\in
Y}\mu\left(\left((K+\sum_{j=1}^p\theta_j (0\times B_j))\cap
(X\times \{y\})\right)+\sum_{i=1}^k\lam_iA_i\right)d\nu(y) =$$
$$\frac{\pt ^p}{\pt \theta_1\dots
\pt\theta_p}\big |_0\int _{y\in Y}
\Theta_X(w)\left(\left(K+\sum_{j=1}^p\theta _j(0\times
B_j)\right)\cap (X\times\{y\})\right)d\nu(y).$$

Let us return now to the case of general $w\in F_X$. To prove the
equality (\ref{diff}) it remains to show that the right hand side
of (\ref{diff}) is continuous with respect to $w\in F_X$ for fixed
$K\in \ck(X\times Y)$ and $B_1,\dots,B_p\in \ck(Y)$. Clearly it is
enough to prove the continuity with respect to $w\in
C^\infty(X\times\PP_+(X^*)^k,|\ome_X|\boxtimes L^{\boxtimes k})$
for any $k=0,1,\dots,n$.

Let us fix a large compact subset $G\subset X$ containing the
projection of $K$ to $X$ in its interior. By Lemma
\ref{techn-lemma} there exist a constant $C$, a compact subset
$\tilde G\subset X$ containing $G$, and an integer $N$ such that
$$w=\sum_{s=1}^\infty \mu_s\otimes h^s_1\otimes \dots\otimes
h^s_k$$ where $\mu_s\in C^\infty(X,|\omega_X|),h^s_i\in
C^\infty(\PP_+(X^*),L)$ and
\begin{eqnarray}
\sum_{s=1}^\infty ||\mu_s||_{C^{k+p+2}(G)}\prod_{i=1}^k
||h^s_i||_{C^{k+p+2}(\PP_+(X^*))}\leq C||w||_{C^N(\tilde G\times
\PP_+(X^*)^k )}.
\end{eqnarray}
Adding to and subtracting from each $h^s_i$ a supporting
functional of the unit Euclidean ball times a constant depending
on $||h^s_i||_{C^2(\PP_+(X^*))}$, we may assume that
$h^s_i=h_{A^s_i}$ is a supporting functional of a convex compact
set $A^s_i$. Thus
$$w=\sum_{s=1}^\infty \mu_s\otimes h_{A^s_1}\otimes \dots\otimes
h_{A^s_k}$$ and
\begin{eqnarray}\label{est1}
\sum_{s=1}^\infty ||\mu_s||_{C^{k+p+2}(G)}\prod_{i=1}^k
||h_{A^s_i}||_{C^{k+p+2}(\PP_+(X^*))}\leq C||w||_{C^N(\tilde
G\times \PP_+(X^*)^k)}.
\end{eqnarray}
Theorem \ref{aff-2-1} implies also that there exist a constant
$C'$, depending on the $B_j$ and $G$, and a compact subset
$G'\subset Y$ such that
\begin{eqnarray}\label{est2}
\big|\big|\frac{\pt ^k}{\pt\lam_1\dots \pt\lam_k}\big |_0
(\mu\boxtimes\nu) \left(K+ \sum_{i=1}^k\lam_i(A_i\times
0)+\sum_{j=1}^p
\theta_j(0\times B_j)\right)\big|\big|_{C^p[0,1]^p}\leq\\
C'||\nu||_{C^{k+p+1}(G')}\cdot
||\mu||_{C^{k+p+1}(G)}\prod_{i=1}^k||h_{A_i}||_{C^{k+p+1}}
\end{eqnarray}
where the function in the left hand side of the inequality is
considered as a function of $(\theta_1,\dots,\theta_p)\in
[0,1]^p$. Hence the function
\begin{eqnarray}\label{eeee}
(\theta_1,\dots,\theta_p)\mapsto \sum_{s=1}^\infty \frac{\pt
^k}{\pt\lam_1\dots \pt\lam_k}\big |_0(\mu_s\boxtimes\nu) \left(K+
\sum_{i=1}^k\lam_i(A_i^s\times 0)+\sum_{j=1}^p \theta_j(0\times
B_j)\right)=\\ \int _{y\in Y}
\Theta_X(w)\left(\left(K+\sum_{j=1}^p\theta _j(0\times
B_j)\right)\cap (X\times\{y\})\right)d\nu(y)
\end{eqnarray}
belongs to $C^p[0,1]^p$, and using (\ref{est1}) and (\ref{est2})
the sum of its $C^p$-norms if the summands in (\ref{eeee})  can be
estimated from above by
\begin{eqnarray*}
C'||\nu||_{C^{k+p+1}(G')}\cdot\sum_s
||\mu_s||_{C^{k+p+1}(G)}\prod_{i=1}^k||h_{A_i^s}||_{C^{k+p+1}}\leq
C'C||\nu||_{C^{k+p+1}(G')}||w||_{C^N(\tilde G)}.
\end{eqnarray*}
Equality (\ref{diff}) is proved, and hence Lemma \ref{mult-1-1}
follows. \qed

For any $\phi\in \sv(X),\, \psi\in \sv(Y)$ we will denote
$M'(\phi,\psi)$ by $\phi\boxtimes \psi$ and call it the exterior
product of $\phi$ and $\psi$. In \cite{alesker-poly} we have
defined the exterior product of polynomial smooth valuations. The
point of this construction is that it extends to smooth valuations
without any assumption of polynomiality.

Let us define now the product on $\sv(V)$. Let
$\Delta:V\hookrightarrow V\times V$ be the diagonal imbedding. For
$\phi,\, \psi\in \sv(V)$ set
$$\phi\cdot \psi:=\Delta^*(\phi\boxtimes \psi)$$
where $\Delta^*$ denotes the restriction of a valuation on
$V\times V$ to the diagonal.
\begin{theorem}\label{mult-1-2}
(1) For $\phi,\, \psi\in \sv(V)$ the product $\phi\cdot \psi \in
\sv(V)$.

(2) The product $SV(V)\times SV(V)\to SV(V)$ is continuous.

(3) Equipped with this multiplication, $\sv(V)$ becomes an
associative commutative unital algebra when the unit is the Euler
characteristic.

(4) The filtration $\{W_\bullet\}$ is compatible with this
multiplication, i.e.
$$W_i\cdot W_j\subset W_{i+j}.$$
\end{theorem}
{\bf Proof.} To prove (1) notice first of all that
$\Delta^*:QV(V\times V)\to QV(V).$ Hence $\sv(V)\cdot \sv
(V)\subset QV(V)$. But since the product commutes with the action
of $Aff(V)$, the product of $Aff(V)$-smooth vectors is
$Aff(V)$-smooth. Hence $\sv(V)\cdot \sv(V)\subset \sv(V)$. The
continuity of the product follows by the same reason.

Let us prove (3) and (4). Using Corollary \ref{aff-3-10} they
reduce to the case of polynomial valuations. But for polynomial
valuations the corresponding statements were proved in
\cite{alesker-poly}. \qed

Let us now describe the associated graded algebra $gr_W \sv(V)$
with respect to the filtration $\{W_i\}$. Remind that $gr_W
\sv(V):=\oplus_{i=0}^n W_i/W_{i+1}$, and it carries the natural
algebra structure.
\begin{theorem}\label{mult-1-3}
The associated graded algebra $gr_W \sv(V)$ is canonically
isomorphic to the graded algebra $C^{\infty}(V, Val^{sm}(V))$ with
the pointwise multiplication on $V$ and the $k$-th graded term of
it is equal to $C^{\infty}(V, Val_k^{sm}(V))$.
\end{theorem}
{\bf Proof.} First let us remind that the isomorphism
$W_k/W_{k+1}$ with $C^{\infty}(V, Val_k(TV))$ is induced by the
map $\Lambda_k\colon W_k\to C^\infty(V,Val_k(TV))$ defined in
Section \ref{aff-3}. Let $\phi\in W_k$. We have
$(\Lambda_k(\phi))(x)(K)= \lim _{r\to +0} r^{-k}\phi(rK+x)$. Thus
the isomorphism of vector spaces follows from Proposition
\ref{aff-3-5}. Now it remains to check that this map is a
homomorphism of algebras. By Corollary \ref{aff-3-10} the result
reduces to the case of polynomial valuations. But for polynomial
valuations the result was proved in \cite{alesker-poly}. \qed
\section{Integration with respect to the normal
cycle.}\label{intcc} In Subsection \ref{summary} we fix some
notation and summarize known relevant facts about construction of
valuations using integration with respect to the normal cycle. The
main new results of this section are contained in Subsection
\ref{main}. These are Theorems \ref{cc-2-1} and \ref{cc-2-2} about
construction of smooth valuations using the integration with
respect to the normal cycle.
\subsection{Main construction and its properties.}\label{summary} \setcounter{subsection}{1}
Let $V$ be a
real vector space of dimension $n$. Then clearly $T^*V=V\times
V^*$. Let $K\in \ck(V)$. Let $x\in K$.
\begin{definition}\label{cc-1-1}
A tangent cone to $K$ at $x$ is a set denoted by $T_xK$ which is
equal to the closure of the set $\{y\in V|\exists \eps>0\, \,
x+\eps y\in K\}$.
\end{definition}
It is easy to see that $T_xK$ is a closed convex cone.
\begin{definition}\label{cc-1-2}
A normal cone to $K$ at $x$ is the set
$$Nor_xK:=\{y\in V^*| \,\, y(x)\geq 0 \,\forall x\in T_xK\}.$$
\end{definition}
Thus $Nor_xK$ is also a closed convex cone.
\begin{definition}\label{cc-1-3}
Let $K\in \ck(V)$. The {\itshape characteristic cycle} of K is the
set $$CC(K):=\cup_{x\in K}Nor_x(K).$$
\end{definition}
\begin{remark}
The notion of the characteristic cycle is not new. First an almost
equivalent notion of normal cycle (see below) was introduced by
Wintgen \cite{wintgen}, and then studied further by Z\"ahle
\cite{zahle} by the tools of geometric measure theory.
Characteristic cycles of subanalytic sets of real analytic
manifolds were introduced by Kashiwara (see
\cite{kashiwara-schapira}, Chapter 9) using the tools of the sheaf
theory, and independently by J. Fu \cite{fu-94} using rather
different tools of geometric measure theory. The elementary
approach described above is sufficient for the purposes of this
article.
\end{remark}

It is easy to see that $CC(K)$ is a closed $n$-dimensional subset
of $T^*V=V\times V^*$ invariant with respect to the multiplication
by non-negative numbers acting on the second factor.
\def\ucc{\underline{CC}}
\def\tcc{\tilde{CC}}
Sometimes we will also use the following notation. Let
$\underline{0}$ denote the zero section of $T^*V$, i.e.
$\underline{0}=V\times\{0\}$. Set
\begin{eqnarray*}
\ucc(K):=CC(K)\backslash\underline{0},\\
\tcc(K):=\ucc/\RR_{>0}.
\end{eqnarray*}
Thus $\tcc(K)\subset \PP_+(T^*V)$. Let us denote by $N(K)$ the
image of $\tcc(K)$ under the involution on $\PP_+(T^*V)$ of the
change of an orientation of a line. $N(K)$ is called the {\itshape
normal cycle} of $K$.

Let us denote by $$p:T^*V\to V$$ the canonical projection. Let us
denote by $o$ the orientation bundle of $V$. Note that a choice of
orientation on $V$ induces canonically an orientation on $CC(K)$
and $N(K)$ for any $K\in \ck(V)$. Let us denote by $\tilde
C^1(T^*V,\Ome^n\otimes p^*o)$ the space of $C^1$-smooth sections
of $\Ome^n\otimes p^*o$ over $T^*V$ such that the restriction of
$p$ to the support of this section is proper.
\begin{theorem}\label{cc-1-4}
For any $\ome\in\tilde C^1(T^*V,\Ome^n\otimes p^*o)$ the map
$\ck(V)\to \CC$ given by $K\mapsto \int_{CC(K)}\ome$ defines a
continuous valuation on $\ck(V)$.
\end{theorem}
The proof of this result can be found in \cite{part3}. However for
a special choice of the form $\ome$ leading to the curvature
measures this theorem was proved much earlier by M. Z\"ahle
\cite{zahle}.

We immediately obtain the following corollary.
\begin{corollary}\label{cc-1-5}
For any $\eta\in C^1(\PP_+(T^*V),\Ome^{n-1}\otimes p^*o)$  the map
$\ck(V)\to \CC$ given by $K\mapsto \int_{N(K)}\eta$ defines a
continuous valuation on $\ck(V)$.
\end{corollary}
We will also need the following statement.
\begin{theorem}[\cite{part3}]\label{cc-1-6}
The map $\ck(V)\times \left(C^1(V,|\ome_V|)\oplus
C^1(\PP_+(V^*),\Ome^{n-1}\otimes p^*o)\right)\to \CC$ given by
$$(K,(\ome,\eta))\mapsto \int_K\ome +\int_{N(K)}\eta$$
is continuous.
\end{theorem}
Theorem \ref{cc-1-6} immediately implies the following corollary.
\begin{corollary}\label{cc-1-7}
(i) The map $C^1(V,|\ome_V|)\oplus
C^1(\PP_+(V^*),\Ome^{n-1}\otimes p^*o)\to CV(V)$ given by
$(\ome,\eta)\mapsto[K\mapsto\int_K\ome +\int_{N(K)}\eta]$ is
continuous.

(ii) For any compact set $G\subset V$ the exists a larger compact
set $\tilde G\subset V$ and a constant $C=C(G)$ such that for any
$(\ome,\eta)\in C^1(V,|\ome_V|)\oplus
C^1(\PP_+(V^*),\Ome^{n-1}\otimes p^*o)$ one has
$$\sup_{K\subset G,K\in\ck(V)}|\int_K\ome+\int_{N(K)}\eta|\leq
C(||\ome||_{C^1(\tilde G)}+||\eta||_{C^1(p^{-1}\tilde G)}).$$
\end{corollary}

\begin{proposition}\label{cc-1-8}
(i) For any $(\ome,\eta)\in C^\infty(V,|\ome_V|)\oplus
C^\infty(\PP_+(V^*),\Ome^{n-1}\otimes p^*o)$ the valuation
$[K\mapsto\int_K\ome+ \int_{N(K)}\eta]$ is smooth, i.e. belongs to
$SV(V)$.

(ii) The induced map $$C^\infty(V,|\ome_V|)\oplus
C^\infty(\PP_+(V^*),\Ome^{n-1}\otimes p^*o)\to SV(V)$$ is
continuous.
\end{proposition}
\def\tx{\tau_{(x,t)}}
{\bf Proof.} Since the construction of integration with respect to
the normal cycle is equivariant with respect to the natural action
of the group $GL(V)$ on all spaces, it is sufficient to prove the
proposition with $SV(V)$ replaced with $QV(V)$ everywhere. For
simplicity we will ignore the summand $\int_K\ome$. The last case
is simpler and it can be considered similarly.

 For any
$(x,t)\in V\times [0,1]$ let us define the map
$\tx:\PP_+(T^*V)\to\PP_+(T^*V)$ by
$$\tx((y,n)):=(ty+x,n).$$ Then we have
$$\int_{N(tK+x)}\eta=\int_{N(K)}\tx^*\eta.$$ Clearly the form
$\tx^*\eta$ depends smoothly on $(x,t)$. This implies part (i) of
the proposition.

Let us prove part (ii). Let us fix a compact set $G\subset V$ and
$N\in \NN$. We have
\begin{eqnarray*}
\sup _{x\in G,K\subset G}||t\mapsto
\int_{N(tK+x)}\eta||_{C^N[0,1]}=\\
\sup _{x\in G,K\subset
G}||t\mapsto\int_{N(K)}\tx^*\eta||_{C^N[0,1]}\leq
C||\eta||_{C^{N+1}(p^{-1}\tilde G)}
\end{eqnarray*}
where $C$ and $\tilde G$ are from Corollary \ref{cc-1-7}(ii). This
implies Proposition \ref{cc-1-8}. \qed

\subsection{Main results}\label{main} \setcounter{theorem}{0}
\def\ctr{C^\infty_{tr}(\PP_+(T^*V),\Ome^{n-1}\otimes p^*o)}
Let us denote by $\ctr$ the subspace consisting of elements of the
space $C^\infty(\PP_+(T^*V),\Ome^{n-1}\otimes p^*o)$ which are
invariant under translations with respect to vectors from $V$.
Elements of this space define translation invariant smooth
valuations.
\begin{theorem}\label{cc-2-1}
Consider the map
$$\CC\cdot vol_V\oplus \ctr\to Val^{sm}(V)$$
given by $(\ome,\eta) \mapsto [K\mapsto \int_K\ome
+\int_{N(K)}\eta].$ This map is onto.
\end{theorem}
{\bf Proof.} Clearly this map commutes with the natural action of
$GL(V)$ on both spaces. It is easy to see that the image of this
map intersect non-trivially each subspace $Val_i^{ev/odd}$ for
$i=0,1,\dots,n$. Hence by Irreducibility Theorem \ref{irr} the
image of this map is dense in $Val^{sm}(V)$. By the
Casselman-Wallach Theorem \ref{casselman-wallach} the image of
this map is closed. Hence it coincides with $Val^{sm}(V)$. \qed

\begin{theorem}\label{cc-2-2}
The map
$$C^\infty(V,|\ome_V|)\oplus
C^\infty(\PP_+({T^*V}),\Ome^{n-1}\otimes p^*o)\to SV(V)$$ is onto.
\end{theorem}

In order to prove this theorem we will introduce a decreasing
filtration on the space $C^\infty(V,|\ome_V|)\oplus
C^\infty(\PP_+{T^*V}),\Ome^{n-1}\otimes p^*o)$ and show that it
maps onto the filtration $W_{\bullet}$ on $SV(V)$. Let us start
with some general considerations.

Let $X$ be a smooth manifold. Let $p:P\to X$ be a smooth bundle.
Let $\Omega^N(P)$ be the vector bundle over $P$ of $N$-forms. Let
us introduce a filtration of $\Omega^N(P)$ by vector subbundles
$W_i(P)$ as follows. For every $y\in P$ set
\begin{multline*} (W_i(P))_y:=\{\omega\in \wedge^NT_y^*P \big|\,
\omega|_F\equiv 0 \mbox{ for all } F\in Gr_N(T_yP)\\
\mbox{ with } \dim(F\cap T_y(p^{-1}p(y)))>N-i\}.\end{multline*}

Clearly we have
$$\Omega^N(P)=W_0(P)\supset W_1(P)\supset \dots \supset W_N(P)\supset
W_{N+1}(P)=0.$$ We will study this filtration in greater detail.

Let us make some elementary observations from linear algebra.
\def\wel{W(L,E)}
Let $L$ be a finite dimensional vector space. Let $E\subset L$ be
a linear subspace. For a non-negative integer $i$ set
$$\wel_i:=\{\omega\in \wedge^NL^*\big|\, \omega|_F\equiv 0\mbox{ for
all } F\subset L \mbox{ with } \dim(F\cap E)>N-i\}.$$ Clearly
$$\wedge^NL^*=\wel_0\supset \wel_1\supset \dots \supset \wel_N\supset
\wel_{N+1}=0.$$
\begin{lemma}\label{filt-1}
There exists canonical isomorphism of vector spaces
$$\wel_i/\wel_{i+1}=\wedge^{N-i}E^*\otimes \wedge^i(L/E)^*.$$
\end{lemma}
{\bf Proof.} First note that for every $0\leq j\leq N$ we have
canonical map
$$\wedge^jE^\perp \otimes \wedge^{N-j}L^*\to \wedge^NL^*$$
given by $x\otimes y\mapsto x\wedge y$. It is easy to see that
$$\wel_i= Im[\oplus_{j\geq i}(\wedge^jE^\perp \otimes
\wedge^{N-j}L^*)\to \wedge^NL^*].$$ Note that the induced map
$$\wedge^iE^\perp \otimes \wedge^{N-i}L^*\to \wel_i/\wel_{i+1}$$
is surjective and factorizes as follows (using the equality
$E^\perp=(L/E)^*$ and the canonical map $L^*\to E^*$)
$$
\Vtriangle<1`1`-1;>[\wedge^i E^\perp \otimes\wedge^{N-i}L^*`
\wel_i/\wel_{i+1}`\wedge^i (L/E)^*\otimes \wedge^{N-i}E^*;``].
$$
Let us check that the obtained map
$$\wedge^i (L/E)^*\otimes \wedge^{N-i}E^*\to \wel_i/\wel_{i+1}$$
is an isomorphism. Let us fix a splitting $L=E\oplus F$. Then
$$\wel_i\simeq \oplus_{j\geq i}\wedge^jF^*\otimes
\wedge^{N-j}E^*.$$ Hence $\wel_i/\wel_{i+1}\simeq
\wedge^iF^*\otimes \wedge^{N-i}E^*\simeq \wedge^i(L/E)^*\otimes
\wedge^{N-i}E^*$. \qed

Let us apply the above construction to the case $P=\PP_+(T^*X)$
with $X$ being a smooth manifold of dimension $n$. The above
construction defines a filtration of the vector bundle
$\Omega^{n-1}(P)$ by vector subbundles. Twisting by the pullback
$p^*o$ of the orientation sheaf $o$ of $X$ we obtain a filtration
$\{W_{\bullet}(P)\}$ by vector subbundles of the vector bundle
$\Omega^{n-1}(P)\otimes p^*o$:
$$\Omega^{n-1}(P)\otimes p^*o=W_0(P)\supset W_0(P)\supset
W_1(P)\supset \dots \supset W_{n-1}(P).$$

Let us denote by $\Omega^{n-1}_{P/X}(P)$ the vector bundle over
$P$ of differential forms along the fibers. (Thus
$\Omega^{n-1}_{P/X}(P)$ is the quotient bundle of
$\Omega^{n-1}(P)$.)
\begin{lemma}\label{filt-2}
For $0\leq i\leq n-1$ there exists a canonical isomorphism
\begin{eqnarray}\label{J}
 J_i\colon W_i(P)/W_{i+1}(P)\tilde \to
\Omega^{n-1-i}_{P/X}\otimes p^*(\wedge^iT^*X)\otimes
p^*o.\end{eqnarray}
\end{lemma}
{\bf Proof.} This is an immediate corollary of Lemma \ref{filt-1}.
\qed

Remind that we denote by $Val^{sm}_i(TX)$ the (infinite
dimensional) vector bundle over $X$ whose fiber over a point $x\in
X$ is equal to the space of translation invariant $i$-homogeneous
$GL_n(\RR)$-smooth valuations on $T_xX$. For any point $x\in X$ we
have the canonical map
\begin{eqnarray}\label{filt-2.5}
C^\infty(\Omega^{n-1-i}(\PP_+(T^*_xX)))\otimes \wedge^iT^*X\otimes
o_{T_x^*X} \to Val_i^{sm}(T_xX)\end{eqnarray}
 where $o_{T_x^*X}$
denotes the orientation sheaf of $T_x^*X$ (this map is given by
integration over a normal cycle). This map induces a continuous
linear map
\begin{eqnarray}\label{Psi}
\Psi_i\colon C^\infty(P,\Omega_{P/X}^{n-1-i}(P)\otimes
p^*(\wedge^iT^*X)\otimes p^*o)\to
C^\infty(X,Val_i^{sm}(TX)).\end{eqnarray}

Let us now apply these constructions to an affine space $V$
(instead of $X$). We will identify $P:=\PP_+(T^*V)\to V$ with
$\PP_+(V^*)\times V$. Then the projection $p\colon
P=\PP_+(V^*)\times V$ is the projection to the second factor.
Consider the following map
$$\Xi\colon C^\infty(V,|\omega_V|)\oplus C^\infty(P,\Omega^{n-1}(P)\otimes p^*o)\to SV(V)$$ which is given by
$$\Xi((\nu,\eta))(K)=\nu(K)+\int_{N(K)}\eta.$$
\begin{proposition}\label{filt-3}
\begin{eqnarray*}
\Xi(C^\infty(V,|\omega_V|))&=&W_n\\
\Xi(C^\infty(P,W_i(P))\oplus C^\infty(V,|\omega_V|))&=& W_i,\,
i=0,1,\dots, n-1.
\end{eqnarray*}
where $W_i$ in the right hand side denotes the $i$-th term of the
filtration on $SV(V)$.
\end{proposition}
{\bf Proof.} The statement is obvious for $i=n$. Assume that
$i<n$. Let us fix for simplicity of notation an orientation on
$V$. Thus the orientation sheaf $o$ becomes trivialized. First let
us show that
\begin{eqnarray}\label{inc}
\Xi(C^\infty(P,W_i(P)))\subset W_i.\end{eqnarray} Fix any
$\omega\in C^\infty(P,W_i(P))$. We have to show that for any $K\in
\ck(V)$ and any $x\in V$
$$\int_{N(tK+x)}\omega =O(t^i) \mbox{ as } t \to +0.$$ This easily
follows from the fact that any such $\omega$ belongs to the space
$$\oplus_{j=i}^{n-1}C^\infty(P,\Omega_{P/V}^{n-1-j}(P))\otimes
\wedge^jV^*.$$ Thus the inclusion (\ref{inc}) is proved. Hence we
obtain a map
\begin{eqnarray}\label{Xi}
\Xi_i\colon C^\infty(P,W_i(P)/W_{i+1}(P))\to W_i/W_{i+1}.
\end{eqnarray}
We will show that $\Xi_i$ is surjective. This will imply
Proposition \ref{filt-3} by the induction in $i$. Remind that by
Proposition \ref{aff-3-5} we have canonical isomorphism
\begin{eqnarray}\label{I}
I_i\colon W_i/W_{i+1} \tilde\to C^\infty (V,Val_i^{sm}(V)).
\end{eqnarray}
One has the following lemma.
\def\aa{C^\infty(P,W_i(P)/W_{i+1}(P))}
\def\bb{W_i/W_{i+1}}
\def\cc{C^\infty(P,\Omega^{n-1-i}_{P/X}(P)\otimes
\wedge^iV^*)}
\def\dd{C^\infty(V,Val_i^{sm}(V))}
\def\ff{\Xi_i}
\def\gg{J_i}
\def\hh{I_i}
\def\kk{\Psi_i}
\begin{lemma}\label{filt-4}
The following diagram is commutative:
\begin{eqnarray}\label{filt-4.5}
\square<1`1`1`1;1200`700>[\aa`\bb`\cc`\dd;\ff`\gg`\hh`\kk]
\end{eqnarray}
 where the maps $J_i,\, \Psi_i,\, \Xi_i,\, I_i$ are defined by
(\ref{J}), (\ref{Psi}), (\ref{Xi}), (\ref{I}) respectively.
\end{lemma}
Let us postpone the proof of Lemma \ref{filt-4} and finish the
proof of Proposition \ref{filt-3}. Since $J_i$ and $I_i$ are
isomorphisms it is enough to prove surjectivity of $\Psi_i$. Let
us consider the Fr\'echet spaces
\begin{eqnarray*}
F_1&:=& C^\infty(\PP_+(V^*),\Omega^{n-1-i}(\PP_+(V^*))\otimes \wedge^i V^*),\\
F_2&:=& Val_i^{sm}(V).
\end{eqnarray*}
Then clearly
\begin{eqnarray*}
C^\infty(P,\Omega_{P/X}^{n-1-i}(P)\otimes \wedge^iV^*)&
=&C^\infty(V,F_1),\\
C^\infty(V,Val_i^{sm})&=& C^\infty(V,F_2).
\end{eqnarray*}
By (\ref{filt-2.5}) we have the canonical map
$$f_i\colon F_1\to F_2.$$
Clearly for any $\psi\in C^\infty(V,F_1)$ and any $y\in V$ one has
$$(\Psi_i\psi)(y)=f_i(\psi(y)).$$
Moreover $f_i$ is surjective by the Irreducibility Theorem
\ref{irr} and the Casselman-Wallach theorem
\ref{casselman-wallach}. Hence $\Psi_i$ is surjective by
Proposition \ref{epi}. \qed

Thus it remains to prove Lemma \ref{filt-4}.

{\bf Proof of Lemma \ref{filt-4}.} Remind that for $\phi\in
W_i/W_{i+1}$ and for all $x\in V,\, K\in \ck(V)$  one has
$$(I_i\phi)(x,K)=\lim_{r\to +0}\frac{1}{r^i}\phi(rK+x).$$
Let us fix $\eta\in C^\infty(P,W_i(P)/W_{i+1}(P))$. Let us fix a
basis $e_1^*,\dots, e_n^*$ in $V^*$. Then we can write
$$J_i(\eta)=\sum_{j_1,\dots,j_i}\eta_{j_1,\dots,j_i}\otimes
e_{j_1}^*\wedge \dots\wedge e_{j_i}^*$$ where
$\eta_{j_1,\dots,j_i}\in C^\infty(P,\Omega^{n-1-i}_{P/X}(P))$.
Then
\begin{eqnarray*}
(I_i(\Xi_i\eta))(x,K)&=&\sum_{j_1,\dots,j_i}\lim_{r\to
+0}\frac{1}{r^i} \int_{N(rK+x)}\eta_{j_1,\dots,j_i}\otimes
e_{j_1}^*\wedge \dots\wedge e_{j_i}^*\\
&=&\sum_{j_1,\dots,j_i}\int_{N(K)}\eta_{j_1,\dots,j_i}|_{p^{-1}(x)}\otimes
e_{j_1}^*\wedge \dots\wedge e_{j_i}^*\\
&=&(\Psi_i(J_i\eta))(x,K).
\end{eqnarray*}
\qed



\begin{thebibliography}{99}
\bibitem{alesker-int}
Alesker, Semyon; Integrals of smooth and analytic functions over
Minkowski's sums of convex sets. MSRI "Convex Geometric Analysis"
34 (1998), 1-15.
\bibitem{alesker-advmath}
Alesker, Semyon; On P. McMullen's conjecture on translation
invariant valuations. Adv. Math. 155 (2000), no. 2, 239--263.
\bibitem{alesker-gafa}
Alesker, Semyon; Description of translation invariant valuations
on convex sets with solution of P. McMullen's conjecture. Geom.
Funct. Anal. 11 (2001), no. 2, 244--272.

\bibitem{alesker-poly}
Alesker, Semyon; The multiplicative structure on polynomial
continuous valuations. Geom. Funct. Anal. 14 (2004), no. 1, 1--26,
also: math.MG/0301148.

\bibitem{part2}
Alesker, Semyon; Theory of valuations on manifolds, II.
math.MG/0503399.

\bibitem{part3}
Alesker, Semyon; Fu, Joseph H.G.; Theory of valuations on
manifolds, III. In preparation.

\bibitem{casselman}
Casselman, William; Canonical extensions of Harish-Chandra modules
to representations of $G$. Canad. J. Math. 41 (1989), no. 3,
385--438.
\bibitem{fu-94}
Fu, Joseph H. G.; Curvature measures of subanalytic sets. Amer. J.
Math. 116 (1994), no. 4, 819--880.
\bibitem{gelfand-vilenkin}
Gelfand, I. M.; Vilenkin, N. Ya.; Generalized functions. Vol. 4.
Applications of harmonic analysis. Translated from the Russian by
Amiel Feinstein. Academic Press [Harcourt Brace Jovanovich,
Publishers], New York-London, 1964 [1977].

\bibitem{godement}
 Godement, Roger; Topologie alg\'ebrique et th\'eorie des faisceaux.
 (French) Troisi\`eme \'edition revue et corrig\'ee.
 Publications de l'Institut de Math\'ematique de l'Universit\'e de Strasbourg, XIII.
 Actualit\'es Scientifiques et Industrielles, No. 1252. Hermann, Paris, 1973.

\bibitem{hadwiger-book}
 Hadwiger, Hugo; Vorlesungen \"uber Inhalt, Oberfl\"ache und Isoperimetrie.
 Springer-Verlag, Berlin-G\"ottingen-Heidelberg 1957.
\bibitem{kashiwara-schapira}
Kashiwara, Masaki; Schapira, Pierre; Sheaves on manifolds. With a
chapter in French by Christian Houzel. Grundlehren der
Mathematischen Wissenschaften [Fundamental Principles of
Mathematical Sciences], 292. Springer-Verlag, Berlin, 1990.
\bibitem{khovanskii-pukhlikov1}
Khovanskii, A. G.; Pukhlikov, A. V.; Finitely additive measures of
virtual polyhedra.
 (Russian) Algebra i Analiz 4 (1992), no. 2, 161--185;
 translation in St. Petersburg Math. J. 4 (1993), no. 2, 337--356.

\bibitem{mcmullen-euler}
McMullen, Peter; Valuations and Euler-type relations on certain
classes of convex polytopes. Proc. London Math. Soc. (3) 35
(1977), no. 1, 113--135.
\bibitem{mcmullen-survey}
McMullen, Peter; Valuations and dissections.
 Handbook of convex geometry, Vol. A, B, 933--988, North-Holland,
 Amsterdam, 1993.
\bibitem{mcmullen-schneider}
 McMullen, Peter; Schneider, Rolf;
Valuations on convex bodies. Convexity and its applications,
170--247, Birkh\"auser, Basel, 1983.

\bibitem{schneider-book}
 Schneider, Rolf; Convex bodies: the Brunn-Minkowski theory.
 Encyclopedia of Mathematics and its Applications, 44.
 Cambridge University Press, Cambridge, 1993.
\bibitem{wallach}
 Wallach, Nolan R.; Real reductive groups. I, II.
  Pure and Applied Mathematics, 132. Academic Press, Inc., Boston, MA,
  1988, 1992.
\bibitem{wintgen}
Wintgen, Peter; Normal cycle and integral curvature for polyhedra
in riemannian manifolds. Differential Geometry (G. Soos and J.
Szenthe, eds.), North-Holland, Amsterdam, 1982.
\bibitem{zahle}
Z\"ahle, Martina; Curvatures and currents for unions of sets with
positive reach. Geom. Dedicata 23 (1987), no. 2, 155--171.
\end{thebibliography}
\end{document}